 \documentclass[11pt]{article}

 \usepackage{mathrsfs,amsfonts,amssymb}
 \usepackage[dvips]{color}

 \newtheorem{definition}{Definition}[section]
 \newtheorem{hypothesis}{Hypothesis}[section]
 \newtheorem{lemma}{Lemma}[section]
 \newtheorem{proposition}{Proposition}[section]
 \newtheorem{theorem}{Theorem}[section]
 \newtheorem{corollary}{Corollary}[section]
 \def\beqlb{\begin{eqnarray}}\def\eeqlb{\end{eqnarray}}
 \def\beqnn{\begin{eqnarray*}}\def\eeqnn{\end{eqnarray*}}
 
 \def\bdefinition{\begin{definition}\sl{}\def\edefinition{\end{definition}}}
 \def\blemma{\begin{lemma}\sl{}\def\elemma{\end{lemma}}}
 \def\bproposition{\begin{proposition}\sl{}\def\eproposition{\end{proposition}}}
 \def\btheorem{\begin{theorem}\sl{}\def\etheorem{\end{theorem}}}
 \def\bcorollary{\begin{corollary}\sl{}\def\ecorollary{\end{corollary}}}

 \def\mbb{\mathbb}\def\mbf{\mathbf}
 \def\mrm{\mathrm}

 \def\proof{\noindent{\it Proof.~~}}
 \def\qed{\hfill$\square$\medskip}
 \def\<{\langle}\def\>{\rangle}

\begin{document}

\vskip1.0cm

\centerline{\Large\textbf{Continuous Time Mixed State Branching Processes}}

\smallskip

\centerline{\Large\textbf{and Stochastic Equations\,\footnote{Supported by the National Key R{\&}D Program of China (No.~2020YFA0712900) and the National Natural Science Foundation of China (No.11531001).}}}

\bigskip

\centerline{ Shukai Chen and Zenghu Li}

\medskip

\centerline{Laboratory of Mathematics and Complex Systems,}

\smallskip

\centerline{School of Mathematical Sciences, Beijing Normal University,}

\smallskip

\centerline{ Beijing 100875, People's Republic of China}

\smallskip

\centerline{E-mails: skchen@mail.bnu.edu.cn and lizh@bnu.edu.cn}

\bigskip

{\narrower{\narrower

\noindent{\it Abstract:} A continuous time mixed state branching process is constructed as the scaling limits of two-type Galton-Watson processes. The process can also be obtained by the pathwise unique solution to a stochastic equation system. From the stochastic equation system we derive the distribution of local jumps and the exponential ergodicity in Wasserstein-type distances of the transition semigroup is given. Meanwhile, we study immigration structures associated with the process and prove the existence of the stationary distribution of the process with immigration.

\smallskip


\smallskip

\par}\par}


\section{Introduction}

\setcounter{equation}{0}

\setcounter{equation}{0}Branching processes were introduced as probabilistic models describing the evolution of populations. The study of branching processes were first initiated by $\mathrm{Bienaym\acute{e}}$ (1845), Galton and Watson (1874), independently, which was known as {\it discrete time and discrete state Galton-Watson processes} (GW-processes). To be more realistic, several naturally generalized processes, {\it continuous time discrete state branching processes} (DB-processes) with or without immigration and {\it continuous time continuous state branching process} (CB-processes) with or without immigration were introduced and studied by researchers.

The DB-processes are continuous time discrete state Markov processes with lifetimes are independent and exponentially distributed random variables. There have been abundant works on the understanding of DB-processes including the construction, properties of moments, limit theorems and so on, we refer to \cite{AN72} for the details. The approach of stochastic equations for branching processes has been developed in the past decades. Let $\mbb{N}=\{0,1,2,...\}$ and $\sharp(\cdot)=\sum_{i}\delta_{i}(\cdot)$ be a counting measure on $\mbb{N}.$ Let $X=\{X_t:t\geq0\}$ be a DB-process with immigration with branching rate $c>0$, offspring distribution $(p_i:i\in\mbb{N})$, immigration rate $\eta>0$ and immigration distribution $(q_i:i\in\mbb{N}).$ The two distributions satisfy $\sum_{k=1}^{\infty}kp_k<\infty$ and $\sum_{k=1}^{\infty}kq_k<\infty$. It is known that $X$ can be obtained as a pathwise unique strong solution of the following stochastic equation:
\beqlb\label{1.1}
X_t=X_0+\int_0^t\int_{\mbb{N}}\int_0^{X_{s-}}(z-1)\,M(\mrm{d}s,\mrm{d}z,\mrm{d}u)
+\int_0^t\int_{\mbb{N}}z\,N(\mrm{d}s,\mrm{d}z),
\eeqlb
where $X_0$ is a random variable taking values in $\mbb{N},$ $M(\mrm{d}s,\mrm{d}z,\mrm{d}u)$ is a Poisson random measure on $(0,\infty)\times\mbb{N}\times(0,\infty)$ with intensity measure $c p_{z}\mrm{d}s\sharp(\mrm{d}z)\mrm{d}u$, $N(\mrm{d}s,\mrm{d}z)$ is a Poisson random measure on $(0,\infty)\times\mbb{N}$ with intensity mesure $\eta q_z\mrm{d}s\sharp(\mrm{d}z)$, and $X_0, M(\mrm{d}s,\mrm{d}z,\mrm{d}u), N(\mrm{d}s,\mrm{d}z)$ are independent of each other. In particular, if $\eta\equiv0,q_z\equiv0$ for all $z\in\mbb{N},$ this reduces to the DB-process. Moreover, here and in the sequel, we understand that for any $b\geq a\geq0$
$$
\int_a^b=\int_{(a,b]}\quad ,\quad \int_a^{\infty}=\int_{(a,\infty)}.
$$

CB-processes were first introduced to applicably model the random evolution of large population dynamics in \cite{J58}. Denote the law on $\mbb{D}([0,\infty),[0,\infty))$ by $\mbb{P}_x$ for each initial value $x\geq0$, the branching property of processes can be described by $\mbb{P}_{x+y}=\mbb{P}_x\ast\mbb{P}_y.$ The semigroup of {\it CB-processes with immigration} (CBI-processes) $(Q_t)_{t\geq0}$ can be characterized uniquely by the Laplace transform:
$$
\int_{[0,\infty)}\mrm{e}^{-\lambda y}\,Q_t(x,\mrm{d}y)=\exp\Big\{-xv_t(\lambda)-\int_0^t\psi(v_s(\lambda))\,\mrm{d}s\Big\},\quad \lambda\geq0,x\geq0,
$$
where for any $\lambda\geq0$, $t\mapsto v_t(\lambda)$ uniquely solves the following equation:
$$
v_t(\lambda)=\lambda-\int_0^t\phi(v_s(\lambda))\,\mrm{d}s,\quad t\geq0,
$$
the branching mechanism $\phi$ and immigration mechanism $\psi$ defined on $[0,\infty)$ take the form of
\beqlb
\phi(z)=az+\alpha z^2+\int_0^{\infty}(\mrm{e}^{-zu}-1+zu)\,m(\mrm{d}u),\qquad
\psi(z)=bz+\int_{0}^{\infty}(1-\mrm{e}^{-zu})\,n(\mrm{d}u)
\eeqlb
with L\'evy measures $m,n$ satisfying $\int_0^{\infty}u\wedge u^2\,m(\mrm{d}u)+\int_0^{\infty}u\,n(\mrm{d}u)<\infty,$ and $a\in\mbb{R}, b,\alpha\geq0.$ In particular, when $\psi\equiv0$, this reduces to the CB-process. There are several ways to construct such processes. \cite{F51} proved that a diffusion process
may arise in a limit theorem of GW-processes. In \cite{KW71}, the authors systematically studied the limit theorems of GW-processes with immigration and characterized the class of the limit process as a CBI-process, the conditions of the main theorem involved iterations of the probability generating functions. Some simpler conditions for the weak convergence were provided in \cite{L06}, and \cite{M09} extended the result of \cite{L06} to a two-type CBI-process. Let $Y=\{Y_t:t\geq0\}$ be a CBI-process, similarly, $Y$ can also be represented as a pathwise unique strong solution to the stochastic equation:
\beqlb\label{1.3}
&&Y_t=Y_0+\int_0^t(b-aY_s)\,\mrm{d}s+\sqrt{2\alpha}\int_0^t\int_0^ {Y_s}\,W(\mrm{d}s,\mrm{d}u)
+\int_0^t\int_0^{\infty}z\,N(\mrm{d}s,\mrm{d}z)\nonumber\\
&&~~~~~~~~~~
+\int_0^t\int_0^{\infty}\int_0^{Y_{s-}}z\,\tilde{M}(\mrm{d}s,\mrm{d}z,\mrm{d}u),
\eeqlb
where $Y_0$ is a random variable taking values in $\mbb{R}_{+}$, $W(\mrm{d}s,\mrm{d}u)$ is a time-space white noise with intensity measure $\mrm{d}s\mrm{d}u,$ $M(\mrm{d}s,\mrm{d}z,\mrm{d}u)$ is a Poisson random measure on $(0,\infty)^3$ with intensity measure $\mrm{d}sm(\mrm{d}z)\mrm{d}u,$ $N(\mrm{d}s,\mrm{d}z)$ is a Poisson random measure on $(0,\infty)^2$ with intensity measure $\mrm{d}sn(\mrm{d}z),$ $\tilde{M}(\mrm{d}s,\mrm{d}z,\mrm{d}u)=M(\mrm{d}s,\mrm{d}z,\mrm{d}u)-\mrm{d}sm(\mrm{d}z)\mrm{d}u$ is the compensated measure of $M(\mrm{d}s,\mrm{d}z,\mrm{d}u)$. Moreover, $Y_0, W, M$ and $N$ are independent of each other. We mention that the moment condition $\int_0^{\infty}z\,n(\mrm{d}z)<\infty$ was removed in \cite{FL10}. The sample paths of $Y$ can also be obtained as a unique strong solution to a stochastic equation driven by Brownian motions and Poisson random measures, one can find the formulation (\ref{1.3}) is more nice to analyse the flows of CBI-processes and other applications, see \cite{DL12} for the specific construction. And we refer to \cite{BL06,DL06,DL12,FL10,L11,L19,P16} for the approach and further properties of the above stochastic equations. Based on stochastic equations established above, \cite{HL16} studied the explicit expression of the distribution of jumps, some applications in Finance can be found in \cite{JMS17}. A two-type CBI-process obtained as a unique strong solution of a stochastic equation system was studied by \cite{M13,M14}.

We can rewrite (\ref{1.3}) without immigration part by extending the $M$ to a Poisson random measure denoted again by $M$ on $(0,\infty)^3\times\mbb{N}$ with intensity $\mrm{d}sm(\mrm{d}z)\mrm{d}u\frac{(\lambda z)^k}{k!}\mrm{e}^{-\lambda z}\sharp(\mrm{d}k)$ for some $\lambda>0$:
\beqlb\label{1.4}
&&Y_t=Y_0-a\int_0^tY_s\,\mrm{d}s+\sqrt{2\alpha}\int_0^t\int_0^ {Y_s}\,W(\mrm{d}s,\mrm{d}u)
+\int_0^t\int_0^{\infty}\int_0^{Y_{s-}}\int_{\mbb{N}}z\,\tilde{M}(\mrm{d}s,\mrm{d}z,\mrm{d}u,\mrm{d}k)\nonumber\\
&&~~~
=Y_0-\phi'(\lambda)\int_0^tY_s\,\mrm{d}s+\sqrt{2\alpha}\int_0^t\int_0^ {Y_s}\,W(\mrm{d}s,\mrm{d}u)\nonumber\\
&&~~~~~~~~~~~~~~
+\int_0^t\int_0^{\infty}\int_0^{Y_{s-}}z\,\tilde{M}^0(\mrm{d}s,\mrm{d}z,\mrm{d}u)
+2\alpha\lambda\int_0^tY_s\,\mrm{d}s\nonumber\\
&&~~~~~~~~~~~~~~
+\int_0^t\int_0^{\infty}\int_0^{Y_{s-}}z\,M^1(\mrm{d}s,\mrm{d}z,\mrm{d}u)
+\int_0^t\int_0^{\infty}\int_0^{Y_{s-}}z\,M^2(\mrm{d}s,\mrm{d}z,\mrm{d}u),
\eeqlb
where $M^0(\mrm{d}s,\mrm{d}z,\mrm{d}u)=M(\mrm{d}s,\mrm{d}z,\mrm{d}u,\{k=0\})$ and
$$
M^1(\mrm{d}s,\mrm{d}z,\mrm{d}u)=M(\mrm{d}s,\mrm{d}z,\mrm{d}u,\{k=1\}),\quad
M^2(\mrm{d}s,\mrm{d}z,\mrm{d}u)=M(\mrm{d}s,\mrm{d}z,\mrm{d}u,\{k\geq2\}).
$$
Recently, \cite{FFA19a} gave another SDE-type description for one-dimensional CB-processes based on (\ref{1.4}) with $a<0, \lambda\geq\lambda^{*}$, here $\lambda^{*}$ is the unique root of $\phi$ on $(0,\infty)$. One of the results in \cite{FFA19a} shows that the last three integrals on the right-hand side of (\ref{1.4}) are identified with the mass that immigrates from the skeleton construction. More precisely, the following stochastic equation system has a unique strong solution:
\beqlb
&&Y_t=Y_0+\int_0^t\Big(2\alpha X_s-\phi'(\lambda)Y_s\Big)\,\mrm{d}s
+\sqrt{2\alpha}\int_0^t\int_0^{Y_s}\,W(\mrm{d}s,\mrm{d}u)\nonumber\\
&&~~~~~~~~~~
+\int_0^t\int_0^{\infty}\int_0^{Y_{s-}}z\,\tilde{M}^0(\mrm{d}s,\mrm{d}z,\mrm{d}u)
+\int_0^t\int_0^{\infty}\int_1^{X_{s-}}z\,M^{3}(\mrm{d}s,\mrm{d}z,\mrm{d}k)\nonumber\\
&&~~~~~~~~~~
+\int_0^t\int_{\mbb{M}}\int_1^{X_{s-}}z_1\,M^{4}(\mrm{d}s,\mrm{d}z,\mrm{d}k),\label{1.5}\\
&&X_t=X_0+\int_0^t\int_{\mbb{M}}\int_1^{X_{s-}}(z_2-1)\,M^{4}(\mrm{d}s,\mrm{d}z,\mrm{d}k),\label{1.6}
\eeqlb
where $\mbb{M}=\mbb{R}_{+}\times\mbb{N}$, $M^3(\mrm{d}s,\mrm{d}z,\mrm{d}k)$ is a Poisson random measure on $(0,\infty)^2\times\{\mbb{N}\setminus\{0\}\}$ with intensity measure $\mrm{d}sz\mrm{e}^{-\lambda z}m(\mrm{d}z)\sharp(\mrm{d}k)$, $M^4(\mrm{d}s,\mrm{d}z,\mrm{d}k)$ is a Poisson random measure on $(0,\infty)\times\mbb{M}\times\{\mbb{N}\setminus\{0\}\}$ with intensity measure $\phi'(\lambda)\mrm{d}s\eta_{z_2}(\mrm{d}z_1)p_{z_2}\sharp(\mrm{d}z_2)\sharp(\mrm{d}k)$, one can see the specific definitions of two distributions $(\eta_{k})_{k\in\mbb{N}}$ and $(p_k)_{k\in\mbb{N}}$ in \cite[p.1127]{FFA19a} and we omitted here. The authors prove that for any $y\geq0$, $\{Y_t:t\geq0\}$ is a weak solution of $(\ref{1.4})$ with initial value $Y_0=y$ if $X_0$ is Poisson distributed with parameter $\lambda y$. Moreover, (\ref{1.5})--(\ref{1.6}) includes the prolific skeleton decomposition when $\lambda=\lambda^{*},$ see \cite{BFM08} for the properties of this special decomposition. And we refer to \cite{FFA19b} for the similar construction of (\ref{1.5})--(\ref{1.6}) in the setting of superprocesses.

Inspired by the formulations (\ref{1.5})--(\ref{1.6}), the first objective of this paper is to construct a two-dimensional branching Markov process $\{(Y_1(t),Y_2(t)):t\geq0\}$ taking values in $\mbb{M}$ obtained as a unique strong solution to a more generalized stochastic equation system than (\ref{1.5})--(\ref{1.6}), and the process is called the {\it continuous time mixed state branching process} (MSB-process). Indeed, the specific form of the stochastic equation system is as follows:
\beqlb
&&Y_1(t)=Y_1(0)-a_{11}\int_0^tY_1(s)\,\mrm{d}s+\int_0^t\sqrt{2\alpha Y_1(s)}\,\mrm{d}B(s)
+\int_0^t\int_0^{Y_1(s-)}\int_\mbb{M}z_1\,\tilde{N}_1(\mrm{d}s,\mrm{d}u,\mrm{d}z)\nonumber\\
&&~~~~~~~~~~~~~~~~~
+a_{21}\int_0^tY_2(s)\,\mrm{d}s+\int_0^t\int_0^{Y_2(s-)}\int_{\mbb{M}_{-1}}z_1\,N_2(\mrm{d}s,\mrm{d}u,\mrm{d}z),\label{1.7}\\
&&Y_2(t)=Y_2(0)+\int_0^t\int_0^{Y_1(s-)}\int_{\mbb{M}}z_2\,N_1(\mrm{d}s,\mrm{d}u,\mrm{d}z)
+\int_0^t\int_0^{Y_2(s-)}\int_{\mbb{M}_{-1}}z_2\,N_2(\mrm{d}s,\mrm{d}u,\mrm{d}z),\nonumber\\
&&\label{1.8}
\eeqlb
where $a_{21}, \alpha\geq0, a_{11}\in\mbb{R}, \mbb{M}_{-1}=\mbb{R}_{+}\times\mbb{N}_{-1}, \mbb{N}_{-1}=\mbb{N}\cup\{-1\},$ $B$ is a standard Brownian motion, $N_1$ and $N_2$ are two Poisson random measures with intensity measures $\mrm{d}sn_1(\mrm{d}z)\mrm{d}u$ and $\mrm{d}sn_2(\mrm{d}z)\mrm{du}$, respectively, $n_1$ and $n_2$ are two L\'evy measures satisfying some moment conditions. Intuitively, there exist interactions between $\{Y_1(t):t\geq0\}$ and $\{Y_2(t):t\geq0\}$,
therefore (\ref{1.7})--(\ref{1.8}) obviously generalize the (\ref{1.5})--(\ref{1.6}). We mention that the Brownian motion $B$ in (\ref{1.7}) can be replaced by a space-time white noise $W$  and the process has the same law for any fixed initial value.

In the literature of the theory of branching processes, the rescaling approach plays a valuable role on
establishing the connection among those branching processes, which leads to the second purpose of this paper and we establish two results. First, for a sequence of GW-processes $\{X_k(n):n\geq0\}_{k\geq1}$ and a positive sequence $\{\gamma_k\}_{k\geq1}$, we show that on proper conditions $\{X_k(\lfloor\gamma_kt\rfloor):t\geq0\}$ converges as $k\rightarrow\infty$ to a DB-process in distribution, where $\lfloor x\rfloor$ denotes the integral part of $x$. Second, for a sequence of two-type GW-processes $\{(Y_{k,1}(n),Y_{k,2}(n)):n\in\mbb{N}\}_{k\geq1}$, we also prove $\{(k^{-1}Y_{k,1}(\lfloor\gamma_kt\rfloor), Y_{k,2}(\lfloor\gamma_kt\rfloor)):t\geq0\}$
converges in distribution to a MSB-process under parallel conditions. The key of two limit theorems above is mainly inspired by \cite{L06,L11,M09}.

The existence of the stationary distribution and ergodic rates are both important topics in the theory of Markov processes. A necessary and sufficient condition for the existence of the stationary distribution of one-dimensional CBI-processes was initiated by \cite{P72}, see also \cite{L11} for a proof. The sufficient condition for the multi-type case can be found in \cite{JKR18+}. The strong Feller property and exponential ergodicity of such processes in the total variation distance were given in
\cite{LM15} by a coupling of CBI-processes constructed by the stochastic equation driven by time-space noises and Poisson random measures, see also \cite{L19}. In a recent work \cite{L20}, the author considered the ergodicities and exponential ergodicities in Wasserstein and total variation distances of Dawson-Watanabe superprocesses with or without immigration, which clearly includes the multi-type CBI-process case. After constructing the MSB-processes, we also want to study the ergodic theory of such processes and we prove the exponential ergodicity in the $L^1$-Wasserstein distance by establishing a upper bound estimates for the variations of the transition probabilities, which is inspired by the similar results of measure-valued branching processes in \cite{L20}. Moreover, by adding the immigration structures, we give a sufficient and necessary condition for the existence of the stationary distribution of {\it MSB-processes with immigration} (MSBI-processes).

The remainder of the paper is organized as follows. In section 2, we prove a weak convergence theorem from GW-processes to DB-processes. In section 3, we obtain the MSB-process arising in a limit theorem of rescaled two-type GW-processes. In section 4, we provide another construction of MSB-processes by stochastic equation systems. The analysis of distributions of jumps is given in section 5. In section 6, we study the estimates for the variations and the exponential ergodicity
both in the $L^1$-Wasserstein distance $W_1$ for the transition semigroup of such processes. And we prove the existence of the stationary distribution of such processes with immigration in section 7.


\section{The construction of DB-processes}
\setcounter{equation}{0}
Let $\{p_j:j\in\mbb{N}\}$ be a probability distribution on $\mbb{N}$, and denote the generating function by $g(z)=\sum_{j=0}^{\infty}p_jz^j$ on $|z|\leq1$. Let $u(z)=a\left(g(z)-z\right)$ for some $a>0$. A Markov process $\{X_t:t\geq0\}$ with state space $\mbb{N}$ is called a DB-process with branching rate $a>0$ and offspring distribution $\{p_j:j\in\mbb{N}\}$ if its transition probabilities $Q_{ij}(t)$ satisfy
\beqlb
\sum_{j=0}^{\infty}Q_{ij}(t)z^j=\Big[\sum_{j=0}^{\infty}Q_{1j}(t)z^j\Big]^i,\quad i\in\mbb{N},\quad t\geq0,\quad z\in[0,1],\label{2.1}
\eeqlb
which implies the branching property of the process. Denote $F(z,t)=\sum_{k=0}^{\infty}Q_{1k}(t)z^k.$ Clearly, $F=(F(\cdot,t):t\geq0)$ satisfies the semigroup property: $F(\cdot,t+s)=F(F(\cdot,t),s)$ for $t,s\geq0$ and is the unique solution of the following differential equation:
\beqlb
\frac{\partial}{\partial t}F(z,t)=u[F(z,t)],\quad F(z,0)=z.\label{2.2}
\eeqlb
We call $F$ the compound semigroup for the DB-process, we refer to \cite[p.106-107]{AN72} for more details.

We now provide a sufficient condition for the weak convergence of GW-processes to the DB-process. Assume that there exists a sequence of GW-processes $\{X_k(n):n\geq0\}_{k\geq1}$ with parameters $\{g_k\}_{k\geq1}$, and let $\{\gamma_k\}_{k\geq1}$ be a sequence of positive numbers. Denote the $n$-step transition probability for $\{X_k(n):n\geq1\}$ by $Q^n_k$, let $\lfloor x\rfloor$ be the integral part of $x.$ One can see that
$$
\sum_{j=0}^{\infty}Q_k^{\lfloor\gamma_kt\rfloor}(i,j)z^j=
[g_k^{\circ\lfloor\gamma_kt\rfloor}(z)]^i:=\big(F_k(z,t)\big)^i,
\quad i\in\mbb{N}, \quad z\in[0,1], \quad t\geq0,
$$
and
\beqnn
&&F_k(z,t)=z+\sum_{i=1}^{\lfloor\gamma_kt\rfloor}(g_k^{\circ i}(z)-g_k^{\circ (i-1)}(z))\\
&&~~~~~~~~~~=z+{\gamma_k}^{-1}\sum_{i=1}^{\lfloor\gamma_kt\rfloor}U_k(F_k(z,\frac{i-1}{\gamma_k}))\\
&&~~~~~~~~~~=z+\int_0^{\frac{\lfloor\gamma_kt\rfloor}{\gamma_k}}U_k(F_k(z,r))\,\mrm{d}r,
\eeqnn
where $g^{\circ n}(z)$ is defined by $g^{\circ n}(z)=g(g^{\circ(n-1)}(z))$ successively with $g^{\circ 0}(z)=z$ and  $U_k(z)=\gamma_k(g_k(z)-z),0\leq z\leq1.$ For the convenience, we formulate the following conditions:

\smallskip

\noindent $(\textbf{A})~\gamma_k\rightarrow\infty$ as $k\rightarrow\infty.$

\smallskip

\noindent $(\textbf{B})$ The sequence $U_k(z)$ is uniformly Lipschitz on [0,1], and converges to a continuous function $u(z)$ as $k\rightarrow\infty.$

\bproposition\label{pr:2.1}
\noindent(i) Suppose that $(\textbf{A,B})$ hold. Then the limit function of sequence $\{U_k(z)\}_{k\geq1}$
has representation $u(z)=a(g(z)-z)$ as $k\rightarrow\infty$ for all $0\leq z\leq1,$ where $a$ is a strictly positive constant, $g(z)$ is a generation function and $ g^{'}(1-)<\infty$.

\noindent(ii) For any given $u(z)=a(g(z)-z),$ there exists a sequence of $\{U_k\}_{k\geq1}$ such that $(\textbf{A,B})$ hold with $U_k(z)\rightarrow u(z).$
\eproposition

\proof
\noindent(i) The desired result is a corollary of Proposition 3.1 (i) later. Indeed, it suffices to consider the offspring distribution  corresponding to two-type GW-processes case satisfying $v_k(\{i,\cdot\})\equiv0$ for all $i\geq1.$

\noindent(ii) For given $g(z)-z=\sum\limits_{i=0}^{\infty}p_iz^i-z$ and $a>0$, set $\gamma_k=ak$ and
$$
p_{ki} =\left\{\begin{array}{ll}
{p_i/k~,}
&\mbox{ if }i\neq1;  \\
&\\
{(p_1-1)/k+1~,} &\mbox{ if }i=1.
\end{array}\right.
$$
Define $U_k(z)=ak(\sum_{i=0}^{\infty}p_{ki}z^i-z)$, it is not hard to see that $U_k$ satisfies condition (\textbf{B}), and converges to $u(z)$ for all $z\in[0,1].$
\qed

\blemma\label{le:2.2}
Suppose that $(\textbf{A,B})$ hold. Then there are constants $\lambda,N\geq0$ such that $ F_k(z,t)\in[z^{\mrm{e}^{\lambda t}},1]$ for every $t\geq0,~z\in[0,1]$ and $k\geq N.$
\elemma

\proof
Let $b_k:=\gamma_k(g^{'}_k(1-)-1)$. Under condition $(\textbf{B})$, there exists $\lambda\geq0$ such that $2|b_k|\leq \lambda$ for all $k\geq1.$ It is not hard to obtain
$$
\sum_{j=1}^{\infty}jQ_k^{\lfloor\gamma_kt\rfloor}(i,j)=ig^{'}_k(1-)^{\lfloor\gamma_kt\rfloor}=i(1+\frac{b_k}{\gamma_k})^{\lfloor\gamma_kt\rfloor}.
$$
Since $\gamma_k\rightarrow\infty$ as $k\rightarrow\infty,$ there is a $N\geq1$ such that for all $k\geq N$
$$
0\leq(1+\frac{b_k}{\gamma_k})^{\frac{\gamma_k}{\lambda}}\leq(1+\frac{\lambda}{2\gamma_k})^{\frac{\gamma_k}{\lambda}}\leq \mrm{e},
$$
then for $t\geq0$ and $k\geq N$
$$
\sum_{j=1}^{\infty}jQ_k^{\lfloor\gamma_kt\rfloor}(i,j)\leq i\exp\bigg\{\frac{\lambda}{\gamma_k}\lfloor\gamma_kt\rfloor\bigg\}
\leq i\mrm{e}^{\lambda t}.
$$
We get the desired result by Jensen's inequality.
\qed

\blemma\label{le:2.3}
Suppose that $(\textbf{A,B})$ hold. For any $c>0,$ we have $F_k(z,t)\rightarrow$ some $F(z,t)$ uniformly on $[0,\mrm{e}^{-c}]\times[0,c]$ as $k\rightarrow\infty$. And the limit function solves (\ref{2.2}).
\elemma

\proof
We may rewrite
\begin{equation}\label{2.3}
F_k(z,t)=z+\int_0^tU_k(F_k(z,r))\,\mrm{d}r-\big(t-\frac{\lfloor\gamma_kt\rfloor}{\gamma_k}\big)U_k
\Big(F_k(z,\frac{\lfloor\gamma_kt\rfloor}{\gamma_k})\Big).
\end{equation}
By Proposition \ref{pr:2.1} and Lemma \ref{le:2.2}, for $\varepsilon\in(0,1]$, we can take $N\geq1$ large enough such that
$$
|U_k(z)-u(z)|\leq\varepsilon,\quad \mrm{e}^{-c\mrm{e}^{\lambda c}}\leq z\leq1,\quad k\geq N.
$$
Denote the last term on the right hand of equation (\ref{2.3}) by $\varepsilon_k(t,z)$, then
$$
|\varepsilon_k(t,z)|\leq\gamma_k^{-1}M,\quad 0\leq t\leq c,\quad 0\leq z\leq \mrm{e}^{-c},
$$
where $M=1+\sup\limits_\Theta|u(z)|, \Theta=\Big\{z:{\mrm{e}^{-c\mrm{e}^{\lambda c}}}\leq z\leq1\Big\}.$
On the other hand, for $n\geq k\geq N$ we set
$$
K_{k,n}(t,z)=\sup\limits_{0\leq r\leq t}\Big|F_n(z,r)-F_k(z,r)\Big|,
$$
it follows that
$$
K_{k,n}(t,z)\leq2(\gamma_k^{-1}M+\varepsilon c)+\rho\int_0^tK_{k,n}(r,z)\mrm{d}r,\quad 0\leq t\leq c,\quad 0\leq z\leq \mrm{e}^{-c},
$$
where $\rho=\sup\limits_\Theta|u^{'}(z)|.$ By applying Gronwall inequality we obtain
$$
K_{k,n}(t,z)\leq2(\gamma_k^{-1}M+\varepsilon c)\mrm{e}^{\rho t},\quad 0\leq t\leq c,\quad 0\leq z\leq \mrm{e}^{-c},
$$
which follows that $F_k(z,t)\rightarrow$ some $F(z,t)$, and the limit function satisfies (\ref{2.2}).
\qed

By (\ref{2.1}), we see that the transition probabilities $Q=\{Q_{ij}(t):i,j\in\mbb{N},t\geq0\}$ of the DB-process $\{X_t:t\geq0\}$ can be determined by
$$
\sum_{j=0}^{\infty}Q_{ij}(t)z^j=(F(z,t))^i,\quad i\geq1, t\geq0.
$$
Based on Proposition \ref{pr:2.1} and Lemma \ref{le:2.3}, by similar arguments of Theorem 2.9 in \cite{L19}, it is not hard to see that the transition probabilities $Q$ of $\{X_t:t\geq0\}$ is a limit of a sequence transition probabilities $\{Q^{\lfloor\gamma_kt\rfloor}_k(i,j):i,j\in\mbb{N},t\geq0\}_{k\geq1}$ associated with GW-processes in the sense of weak convergence under conditions $(\textbf{A,B})$, which indeed implies another construction of DB-processes by rescaling approach.

\section{The construction of MSB-processes}
\setcounter{equation}{0}
For a two-type GW-process $\{Y(n)=(Y_1(n),Y_2(n)):n\in\mbb{N}\},$  we define two corresponding generation functions
$$
g_j(s_1,s_2)=\sum_{i\in{\mbb{N}}^2}p_j(i)s_1^{i_1}s_2^{i_2},\quad p_j(i)=\mbf{P}(Y(1)=i|Y(0)=\mbf{e}_j), \quad j=1,2
$$
for $i=(i_1,i_2)\in\mbb{N}^2$ and $s_1,s_2\in[0,1],$ where $\mbf{e}_1=(1,0),\mbf{e}_2=(0,1).$ Note that $g_j^{\circ n}=g_j(g_1^{\circ(n-1)},g_2^{\circ(n-1)})$ for $j=1,2$ and $n\geq1$. It is known that the one-step transition matrix $P(i,j)$ of $\{Y(n):n\geq0\}$ is uniquely determined by
$$
\sum_{j\in{\mbb{N}}^2}P(i,j)s_1^{j_1}s_2^{j_2}=\Big[g_1(s_1,s_2)\Big]^{i_1}\Big[g_2(s_1,s_2)\Big]^{i_2}.
$$
Let $\{Y_k(n)= (Y_{k,1}(n),Y_{k,2}(n)): n\in\mbb{N}\}$ be a sequence of two-type GW-processes corresponding to the $G_k=(g_{k,1},g_{k,2}),$ and $\{\gamma_k\}$ be a sequence of positive numbers.
Denote the $n-$step transition probability for $\{(k^{-1}Y_{k,1}(n),Y_{k,2}(n)):n\geq0\}$ by $P_k^n,$ it is not hard to see that for $t\geq0$ and $\lambda\in\mbb{R}_{+}^2$
$$
\int_{\mbb{M}_k}\mrm{e}^{-\langle\lambda,y\rangle}P_k^{\lfloor\gamma_k t\rfloor}(x,\mrm{d}y)=\Big[g_{k,1}^{\circ\lfloor\gamma_k t\rfloor}(\mrm{e}^{-\lambda_1/k},\mrm{e}^{-\lambda_2})\Big]^{kx_1}\Big[g_{k,2}^{\circ\lfloor\gamma_k t\rfloor}(\mrm{e}^{-\lambda_1/k},\mrm{e}^{-\lambda_2})\Big]^{x_2},\quad x\in\mbb{M}_k,
$$
where $\mbb{M}_k:=\{(i/k,j):(i,j)\in\mbb{N}^2\}.$ For $\lambda=(\lambda_1,\lambda_2)\in\mbb{R}_{+}^2$, we define a vector function $V_k(t,\lambda)=(V_{k,1}(t,\lambda),V_{k,2}(t,\lambda))$ with
\beqlb
V_{k,1}(t,\lambda)=-k\log g_{k,1}^{\circ \lfloor\gamma_k t\rfloor}(\mrm{e}^{-\lambda_1/k},\mrm{e}^{-\lambda_2}),\quad
V_{k,2}(t,\lambda)=-\log g_{k,2}^{\circ \lfloor\gamma_k t\rfloor}(\mrm{e}^{-\lambda_1/k},\mrm{e}^{-\lambda_2}),\label{3.1}
\eeqlb
then we can rewrite
$$
\int_{\mbb{M}_k}\mrm{e}^{-\langle\lambda,y\rangle}P_k^{\lfloor\gamma_k t\rfloor}(x,\mrm{d}y)=\exp\{-\langle x,V_{k}(t,\lambda)\rangle\},\quad x\in\mbb{M}_k
$$
and $V_k(t,\lambda)$ satisfies the following equation
\beqlb
V_{k,1}(t,\lambda)=\lambda_1-\int_0^{\frac{\lfloor\gamma_kt\rfloor}{\gamma_k}}\bar{\Phi}_{k,1}(V_k(s,\lambda))\,\mrm{d}s,\quad
V_{k,2}(t,\lambda)=\lambda_2-\int_0^{\frac{\lfloor\gamma_kt\rfloor}{\gamma_k}}\bar{\Phi}_{k,2}(V_k(s,\lambda))\,\mrm{d}s,\label{3.2}
\eeqlb
where two functions $\bar{\Phi}_{k,1}$ and $\bar{\Phi}_{k,2}$ are defined on $\mbb{R}_{+}^2$ taking the form of
\beqnn
&&\bar{\Phi}_{k,1}(\lambda_1,\lambda_2)=k\gamma_k\log[1-(k\gamma_k)^{-1}\Phi_{k,1}(\lambda_1,\lambda_2)\mrm{e}^{\frac{\lambda_1}{k}}],\\
&&\Phi_{k,1}(\lambda_1,\lambda_2)=k\gamma_k[\mrm{e}^{-\frac{\lambda_1}{k}}-g_{k,1}(\mrm{e}^{-\frac{\lambda_1}{k}},\mrm{e}^{-\lambda_2})],\\
&&\bar{\Phi}_{k,2}(\lambda_1,\lambda_2)=\gamma_k\log[1-\gamma_k^{-1}\Phi_{k,2}(\lambda_1,\lambda_2)\mrm{e}^{\lambda_2}],\\
&&\Phi_{k,2}(\lambda_1,\lambda_2)=\gamma_k[\mrm{e}^{-\lambda_2}-g_{k,2}(\mrm{e}^{-\frac{\lambda_1}{k}},\mrm{e}^{-\lambda_2})].
\eeqnn
We further define two functions $\Phi_1$ and $\Phi_2$ on $[0,\infty)^2$:
\beqlb
&&\Phi_1(\lambda_1,\lambda_2)=-a_{11}\lambda_1-\alpha\lambda_1^2- \int_{\mbb{M}}(\mrm{e}^{-\langle\lambda,z\rangle}-1+\lambda_1z_1)\,n_1(\mrm{d}z),\label{3.3}\\
&&\Phi_2(\lambda_1,\lambda_2)=a_{21}\lambda_1+\int_{\mbb{M}_{-1}}(1-\mrm{e}^{-\langle\lambda,z\rangle})\,n_2(\mrm{d}z),\label{3.4}
\eeqlb
where $a_{11}$ is a constant, $a_{21},\alpha\geq0$, $n_1$ and $n_2$ are $\sigma$ finite measures on $\mbb{M}$ and ${\mbb{M}_{-1}}$,  supported by
$ \mbb{M}\setminus\{\mbf{0}\}$ and ${\mbb{M}_{-1}}\setminus\{\mbf{0}\},$ respectively, such that $n_2\Big(\mbb{R}_{+}\times\{-1\}\Big)<\infty$ and
$$
\int_{\mbb{M}}(z_1\wedge z_1^2+z_2)\,n_1(\mrm{d}z)+\int_{\mbb{M}_{-1}}(z_1+z_2)\,n_2(\mrm{d}z)<\infty.
$$
For convenience, let us consider the following conditions:

\noindent $(\textbf{A})~\gamma_k\rightarrow\infty.$

\smallskip

\noindent $(\textbf{C})$ The sequence $\{\Phi_{k,1}(\lambda_1,\lambda_2)\}_{k\geq1}$ is uniformly Lipschitz in $(\lambda_1,\lambda_2)$ on each bounded rectangle, and converges to a continuous function as $k\rightarrow\infty$.

\smallskip

\noindent $(\textbf{D})$ The sequence $\{\mrm{e}^{\lambda_2}\Phi_{k,2}(\lambda_1,\lambda_2)\}_{k\geq1}$ is uniformly Lipschitz in $(\lambda_1,\lambda_2)$ on each bounded rectangle, and converges to a continuous function as $k\rightarrow\infty$.

\smallskip

\bproposition\label{pr:3.1}

\noindent(i) Assume that $(\textbf{A,C,D})$ hold. Then the limit function $\Phi_1(\lambda_1,\lambda_2)$ has the representation (\ref{3.3}), and the limit $\Phi_2(\lambda_1,\lambda_2)$ of has the representation (\ref{3.4});

\noindent(ii) For any $\Phi_1$ and $\Phi_2$ given by (\ref{3.3}) and (\ref{3.4}), there are sequences $\{\gamma_k\}$ and $\{(g_{k,1},g_{k,2})\}$ as above such that $(\textbf{A,C,D})$ hold with $\Phi_{k,1}(\lambda_1,\lambda_2)\rightarrow\Phi_1(\lambda_1,\lambda_2)$ and $\mrm{e}^{\lambda_2}\Phi_{k,2}(\lambda_1,\lambda_2)
\rightarrow\Phi_2(\lambda_1,\lambda_2)$ for all $(\lambda_1,\lambda_2)\in\mbb{R}_{+}^2$ as $k\rightarrow\infty.$
\eproposition

\proof
${\bf(i)}$ We first prove the representation result for $\Phi_2$, which is inspired by \cite{M09}.

\noindent{\bf (1)} Fix $k\geq1.$ Let $\mbb{M}_{-1,k}=\{(i/k,j-1):\{i,j\}\in\mbb{N}^2\}$, $\rho_k$ be the measure defined by:
$$
\rho_k(\cdot)=\gamma_k\sum_{i,j=0}^{\infty}v_k(\{i,j\})\delta_{(i/k,j-1)}(\cdot),
$$
where $v_k$ be the probability measure on $\mbb{N}^2$ corresponding to $g_{k,2},$ then $\rho_k$ is a finite measure on $\mbb{M}_{-1,k}$. Let $l(z)=(z_1+|z_2|)\wedge1,\varrho_k=\int_{\mbb{M}_{-1}}l(z)\,\rho_k(\mrm{d}z)$. If $\varrho_k>0,$ define $ P_k(\mrm{d}z)=(l(z)/\varrho_k)\rho_k(\mrm{d}z),$
if $\varrho_k=0,$ we let $P_k(\cdot)$ be the Dirac measure at some point $z_0\in \mbb{M}_{-1,k}\setminus\{0\}$. In both cases,we have that $P_k(\cdot)$ is a probability measure on $\mbb{M}_{-1}$ and
\beqlb\label{3.5}
\mrm{e}^{\lambda_2}\Phi_{k,2}(\lambda_1,\lambda_2)=\varrho_k \int_{{\mbb{M}_{-1}}\setminus\{0\}}
(1-\mrm{e}^{-\langle\lambda,z\rangle})l(z)^{-1}\,P_k(\mrm{d}z).
\eeqlb
{\bf (2)} Let $\mbb{M}_{-1}^{\Delta}=\mbb{M}_{-1}\cup\{\Delta\}$ be the one-point compactification of $\mbb{M}_{-1}$. Then $\{P_k\}$ is the sequence of probability measures on $\mbb{M}_{-1}^{\Delta}$, so it is relatively compact. Choose any subsequence denoted again by $\{P_k\},$ which converges to a probability $P$ on $\mbb{M}_{-1}^{\Delta}$. Let
$ E=\{\varepsilon\mid P(\|z\|=\varepsilon)=0\},$ for $\varepsilon\in E,$ we define a compact space of $\mbb{M}_{-1}^{\Delta}$ by $Q:=\{z\in \mbb{M}_{-1}^{\Delta},\|z\|\leq\varepsilon\}.$ We have
\beqlb\label{3.6}
\mrm{e}^{\lambda_2}\Phi_{k,2}(\lambda_1,\lambda_2)=\varrho_ka^{(1)}_{k,\varepsilon}\lambda_1+\varrho_ka^{(2)}_{k,\varepsilon}\lambda_2
-\varrho_k(I_{k,\varepsilon}+J_{k,\varepsilon}),
\eeqlb
where
\beqnn
&&a_{k,\varepsilon}^{(i)}=\int_{Q\setminus\{0\}}\chi(z_i)l(z)^{-1}\,P_k(\mrm{d}z),\quad i=1,2,\quad \chi(z_i)=(1\wedge z_i)\vee(-1),\\
&&J_{k,\varepsilon}=\int_{\mbb{M}_{-1}\setminus Q}(\mrm{e}^{-\langle\lambda,z\rangle}-1)l(z)^{-1}\,P_k(\mrm{d}z),\\
&&I_{k,\varepsilon}=\int_{Q\setminus\{0\}}\Big(\mrm{e}^{-\langle\lambda,z\rangle}-1+\chi(z_1)\lambda_1+\chi(z_2)\lambda_2\Big)l(z)^{-1}\,P_k(\mrm{d}z).
\eeqnn
It is not hard to see that $a_{k,\varepsilon}^{(2)}=0, \lim\limits_{E\ni\varepsilon\downarrow0} \lim\limits_{k\rightarrow\infty}I_{k,\varepsilon}=0.$

\noindent{\bf (3)} Fix $\varepsilon\in E,0<\varepsilon<1.$ If $\lim\inf\limits_{k\rightarrow\infty}\varrho_k=0,$ then $\Phi_2=0.$ If $ \lim\inf\limits_{k\rightarrow\infty}\varrho_k>0,$ then there exists a subsequence denoted again by $\{\varrho_k\}$ converges to $\varrho\in(0,\infty].$ We can prove that $P(\{\Delta\})=0.$ Actually,
$$
\int_{\mbb{M}_{-1}\setminus Q}\mrm{e}^{-\langle\lambda,z\rangle}l(z)^{-1}\,P_k(\mrm{d}z)\rightarrow~~\mrm{some}\quad L(\lambda),
$$
where $L(\lambda)$ is a continuous function. Since $l(z)^{-1} \rightarrow 1$ as $z\rightarrow\Delta$ and is the continuous function of $z\in {\mbb{M}_{-1}^{\Delta}}\setminus Q$, it follows from monotone convergence theorem that
\beqlb
&&\int_{{\mbb{M}_{-1}}\setminus Q}l(z)^{-1}\,P(\mrm{d}z)=\lim_{k\rightarrow\infty}\int_{{\mbb{M}_{-1}}\setminus\{0\}}l(z)^{-1}\,P_k(\mrm{d}z)=L(0),\nonumber\\
&&\int_{{\mbb{M}_{-1}^{\Delta}}\setminus Q}l(z)^{-1}\,P(\mrm{d}z)=\lim_{n\rightarrow\infty}\int_{{\mbb{M}_{-1}}\setminus Q}
\mrm{e}^{-\langle\frac{1}{n},z\rangle}l(z)^{-1}\,P(\mrm{d}z)=\lim_{n\rightarrow\infty}L(\frac{1}{n})=L(0).\nonumber
\eeqlb
Define $v(\{0\})=0,$ and $v(\mrm{d}z)=l(z)^{-1}P(\mrm{d}z)$ on $\{z\in {\mbb{M}_{-1}}:\|z\|>0\}$, we have
$$
\lim\limits_{E\ni\varepsilon\downarrow0}\lim\limits_{k\rightarrow\infty}J_{k,\varepsilon}=\int_{\mbb{M}_{-1}}(\mrm{e}^{-\langle\lambda,z\rangle}-1)\,v(\mrm{d}z),
\quad \lim\limits_{E\ni\varepsilon\downarrow0}\lim\limits_{k\rightarrow\infty}a_{k,\varepsilon}^{(1)}:=a^{(1)}\geq0.
$$
Based on calculations above, we obtain that
\beqlb\label{3.7}
\frac{1}{\varrho}\Phi_2(\lambda_1,\lambda_2)=a^{(1)}\lambda_1+\int_{\mbb{M}_{-1}}(1-\mrm{e}^{-\langle\lambda,z\rangle})v(\mrm{d}z).
\eeqlb
{\bf (4)} Now we need to verify that $1/\varrho>0.$ If not, the r.h.s of (\ref{3.7}) is the zero function. Assume $ v(\mbb{M}\setminus\{0\})=0,$ then $P(\{0\})=1,$ and by the representations of $a_{k,\varepsilon}^{(i)},$ we can have $ a^{(1)}=\lim\limits_{E\ni\varepsilon\downarrow0}
P(\|z\|\leq\varepsilon)=1,$ but it follows from the representation of $ \textit{L}\acute{e}\textit{vy-Khintchine}$ type functions that the parameters are unique, therefore $1/\varrho=0$ is impossible. Let $n_2(\cdot)=\varrho v(\cdot),a_{21}=\varrho a^{(1)},$ then we can rewrite $\Phi_2$ as follows:
\beqlb
\Phi_2(\lambda_1,\lambda_2)=a_{21}\lambda_1+\int_{\mbb{M}_{-1}}(1-\mrm{e}^{-\langle\lambda,z\rangle})\,n_2(\mrm{d}z).\nonumber
\eeqlb
By applying the monotone convergence theorem, we get
\beqlb
\frac{\partial\Phi_2}{\partial\lambda_1}(0+,0)=a_{21}+\int_{\mbb{M}_{-1}}z_1\,n_2(\mrm{d}z),\quad \frac{\partial\Phi_2}{\partial\lambda_2}(0,0+)=\int_{\mbb{M}_{-1}}z_2\,n_2(\mrm{d}z).\nonumber
\eeqlb
From the fact that $\Phi_2$ is locally Lipschitz, we have $\int_{\mbb{M}_{-1}}(z_1+z_2)\, n_2(\mrm{d}z)<\infty$ and $n_2(\mbb{R}_{+}\times\{-1\})<\infty.$

Now we treat with $\Phi_{k,1}.$ Fix $k\geq1,$ let $W=[-1,\infty)\times\mbb{N},
W_k=\{((i-1)/k,j):(i,j)\in\mbb{N}^2\}.$ We define $\tilde{\rho}_k$ by:
$$
\tilde{\rho}_k(\cdot)=k\gamma_k\sum\limits_{i,j=0}^{\infty}\tilde{v}_k(\{i,j\})\delta_{(\frac{i-1}{k},j)}(\cdot),\nonumber
$$
where $\tilde{v}_k$ is the probability measure on $\mbb{N}^2$ corresponding to $g_{k,1}$. Let $\tilde{l}(z)=(z_1^2+z_2)\wedge1,$ we have
\beqnn
\mrm{e}^{\lambda_1/k}\Phi_{k,1}(\lambda_1,\lambda_2)=\tilde{\beta}_{k,1}\lambda_1+\tilde{\varrho}_k\tilde{a}_{k,\varepsilon}^{(2)}\lambda_2-
\tilde{\varrho}_k\Big(\tilde{a}_{k,\varepsilon}^{(1)}\lambda_1^2+\tilde{I}_{k,\varepsilon}+\tilde{J}_{k,\varepsilon}\Big),
\eeqnn
where
\beqnn
&&\tilde{\beta}_{k,1}=\int_W\chi(z_1)\,\tilde{\rho}_k(\mrm{d}z),\quad \tilde{a}_{k,\varepsilon}^{(1)}=\frac{1}{2}\int_{Q\setminus\{0\}}
\chi(z_1)^2\tilde{l}(z)^{-1}\,\tilde{P}_k(\mrm{d}z),\\
&&\tilde{a}_{k,\varepsilon}^{(2)}=\int_{Q\setminus\{0\}}\chi(u_2)\tilde{l}(u)^{-1}\,\tilde{P}_k(\mrm{d}u),\quad h(u,\lambda)=\mrm{e}^{-\langle\lambda,z\rangle}-1+\chi(z_1)\lambda_1,\\
&&\tilde{J}_{k,\varepsilon}=\int_{W\setminus Q}h(z,\lambda)\tilde{l}(z)^{-1}\,\tilde{P}_k(\mrm{d}z),\\
&&\tilde{I}_{k,\varepsilon}=\int_{Q\setminus\{0\}}\Big(h(z,\lambda)+\chi(z_2)\lambda_2-\frac{1}{2}\chi(z_1)^2\lambda_1^2\Big)\tilde{l}(z)^{-1}\,\tilde{P}_k(\mrm{d}z).
\eeqnn
Note that $\tilde{a}_{k,\varepsilon}^{(2)}=0.$ By a similar argument in \cite{M09}, we have
\beqnn
\Phi_1(\lambda_1,\lambda_2)=-a_{11}\lambda_1-\alpha\lambda_1^2-\int_{\mbb{M}}(\mrm{e}^{-\langle\lambda,z\rangle}-1+\lambda_1z_1)\,n_1(\mrm{d}z),
\eeqnn
and $n_1$ satisfies
$$
\int_{\mbb{M}}\big(z_1\wedge z_1^2+z_2\big)\,n_1(\mrm{d}z)<\infty.
$$

\noindent${\bf(ii)}$ Given function by (\ref{3.4}), we set $\bar{D}_k=\{z\in {\mbb{M}_{-1}}:z_1>\frac{1}{\sqrt{k}}\},$ $\bar{\gamma}_{k,1}=a_{21}k,
\bar{\gamma}_{k,2}=n_2(\bar{D}_k),$ $ \bar{\mrm{g}}_{k,1}(x_1,x_2)=x_1x_2$ and $\bar{\mrm{g}}_{k,2}(x_1,x_2)=\bar{\gamma}_{k,2}^{-1}
\int_{\bar{D}_k}x_1^{kz_1}x_2^{(z_2+1)}\, n_2(\mrm{d}z).$ Then if we set sequences $\bar{ \gamma}_k=\bar{\gamma}_{k,1}+\bar{\gamma}_{k,2}$ and $\bar{\mrm{g}}^{(2)}_k=\bar{\gamma}_k^{-1}(\bar{\gamma}_{k,1}\bar{\mrm{g}}_{k,1}+\bar{\gamma}_{k,2}\bar{\mrm{g}}_{k,2})$,
we find the sequences $\{\bar{\gamma}_k\}$ and $\{\bar{\mrm{g}}_k^{(2)}\}$ such that conditions $(\textbf{A,D})$ hold with $\mrm{e}^{\lambda_2}\Phi_{k,2}(\lambda_1,\lambda_2)\rightarrow\Phi_2(\lambda_1,\lambda_2)$ for all $(\lambda_1,\lambda_2)\in\mbb{R}_{+}^2$ as $k\rightarrow\infty.$

On the other hand, given function by (\ref{3.3}). Firstly, we set $\tilde{\gamma}_{k,1}=a_{11},\tilde{\mrm{g}}_{k,1}(x_1,x_2)=1.$ Secondly, set $\tilde{\gamma}_{k,2}=(2\alpha+1)k,\tilde{\mrm{g}}_{k,2}(x_1,x_2)=x_1+\frac{\alpha}{2\alpha+1}(1-x_1)^2.$
Thirdly, let $\tilde{D}_k=\{z\in\mbb{M}:z_1>\frac{1}{\sqrt{k}},z_2>\frac{1}{\sqrt{k}}\},$ and $\sigma_k=\int_{\tilde{D}_k}(z_1-\frac{1}{k})\,
n_1(\mrm{d}z),$ set the sequences
\beqnn
&&\tilde{\gamma}_{k,3}=\sigma_k+\frac{1}{k}n_1(\tilde{D}_k)+1,\\
&&\tilde{\mrm{g}}_{k,3}(x_1,x_2)=\frac{1}{k\tilde{\gamma}_{k,3}}\int_{\tilde{D}_k}x_1^{kz_1}x_2^{z_2}\,n_1(\mrm{d}z)+\frac{\sigma_k+1}{\tilde{\gamma}_{k,3}}
\Big(x_1+\frac{\sigma_k}{\sigma_k+1}(1-x_1)\Big).
\eeqnn
Fourthly, if we let $\tilde{\gamma}_k=\tilde{\gamma}_{k,1}+\tilde{\gamma}_{k,2}+\tilde{\gamma}_{k,3}$ and $\tilde{\mrm{g}}_k^{(1)}=\tilde{\gamma}_k^{-1}(\tilde{\gamma}_{k,1}\tilde{\mrm{g}}_{k,1}+\tilde{\gamma}_{k,2}\tilde{\mrm{g}}_{k,2}+\tilde{\gamma}_{k,3}\tilde{\mrm{g}}_{k,3}),$
we find that the sequences $\{\tilde{\gamma}_k\}$ and $\{\tilde{\mrm{g}}_k^{(1)}\}$ such that conditions $(\textbf{A,C})$ hold with $\Phi_{k,1}(\lambda_1,\lambda_2)\rightarrow\Phi_1(\lambda_1,\lambda_2)$ for all $(\lambda_1,\lambda_2)\in\mbb{R}_{+}^2$ as $k\rightarrow\infty.$

Now let $\gamma_k=\tilde{\gamma}_k+\bar{\gamma}_k$ and
\beqnn
&&\mrm{g}_{k,1}(x_1,x_2)=\gamma_k^{-1}(\tilde{\gamma}_k\tilde{\mrm{g}}_k^{(1)}(x_1,x_2)+\bar{\gamma}_kx_1),\\ &&\mrm{g}_{k,2}(x_1,x_2)=\gamma_k^{-1}(\bar{\gamma}_k\bar{\mrm{g}}_k^{(2)}(x_1,x_2)+\tilde{\gamma}_kx_2),
\eeqnn
we finally find the common sequences $\{\gamma_k\}$ and $\{(\mrm{g}_{k,1},\mrm{g}_{k,2})\}$ such that conditions $(\textbf{A,C,D})$ hold with $\Phi_{k,1}(\lambda_1,\lambda_2)\rightarrow\Phi_1(\lambda_1,\lambda_2)$ and $\mrm{e}^{\lambda_2}\Phi_{k,2}(\lambda_1,\lambda_2)\rightarrow\Phi_2(\lambda_1,\lambda_2)$ for all $(\lambda_1,\lambda_2)\in\mbb{R}_{+}^2$ as $k\rightarrow\infty.$
\qed

\bproposition\label{pr:3.2}
Assume that $(\textbf{A,C,D})$ hold. Then for any $a\geq0$ we have $V_k(t,\lambda)\rightarrow$ some $V(t,\lambda)$ uniformly on $[0,a]^3$ as $k\rightarrow\infty$, and the limit function solves the following integral equations:
\beqlb
V_1(t,\lambda)=\lambda_1+\int_0^t\Phi_1(V(s,\lambda))\,\mrm{d}s,\quad
V_2(t,\lambda)=\lambda_2+\int_0^t\Phi_2(V(s,\lambda))\,\mrm{d}s.\label{3.8}
\eeqlb
Moreover, suppose that $(\Phi_1,\Phi_2)$ is given by (\ref{3.3})--(\ref{3.4}). Then for any $\lambda\in\mbb{R}_{+}^2,$ the solution $t\mapsto V(t,\lambda)$ of (\ref{3.8}) is unique, and the solution satisfies the semigroup property:
\beqlb
V(r+t,\lambda)=V(r,V(t,\lambda)),\quad r,t\geq0.\label{3.9}
\eeqlb
\eproposition

\proof  By a similar argument of Lemma 2.2 in \cite{L19}, it follows from Proposition \ref{pr:3.1} that $\bar{\Phi}_{k,1}\rightarrow-\Phi_1$ and $\bar{\Phi}_{k,2}\rightarrow-\Phi_2$ uniformly on each bounded rectangle, respectively, as $k\rightarrow\infty$. We can rewrite
\beqnn
&&V_{k,1}(t,\lambda)=\lambda_1-\int_0^t\bar{\Phi}_{k,1}(V_k(s,\lambda))\,\mrm{d}s+\varepsilon_{k,1}(t,\lambda),\\ &&V_{k,2}(t,\lambda)=\lambda_2-\int_0^t\bar{\Phi}_{k,2}(V_k(s,\lambda))\,\mrm{d}s+\varepsilon_{k,2}(t,\lambda),
\eeqnn
where
$$
\varepsilon_{k,i}(t,\lambda)=-(t-\gamma_k^{-1}\lfloor\gamma_kt\rfloor)\bar{\Phi}_{k,i}(V_k(\gamma_k^{-1}\lfloor\gamma_kt\rfloor,\lambda)),\quad i=1,2.
$$
For $\theta\in\mbb{R}_{+}^2,$ it is not hard to obtain that
\beqnn
&&\int_{\mbb{M}_k}(\theta_1y_1+\theta_2y_2)\,P^{\lfloor\gamma_kt\rfloor}_k(x,\mrm{d}y)\\
&&=\Big(k^{-1}\theta_1g'_{k,1,(1)}(1)+\theta_2g'_{k,1,(2)}(1), k^{-1}\theta_1g'_{k,2,(1)}(1)+\theta_2g'_{k,2,(2)}(1)\Big)
A_k^{(\lfloor\gamma_kt\rfloor-1)}(kx_1,x_2)^\top\\
&&=\bigg(k^{-1}\Big[\theta_1(\gamma^{-1}_k\Phi'_{k,1,(1)}(0)+1)+\theta_2\gamma_k^{-1}\Phi'_{k,1,(2)}(0)\Big],\,\,\theta_1\gamma^{-1}_k\Phi'_{k,2,(1)}(0)\\
&&~~~~~~~~~~~~~~~~~~~~~~~~~~~~~~~~
+\theta_2(\gamma_k^{-1}\Phi'_{k,1,(2)}(0)+1)\bigg)A_k^{(\lfloor\gamma_kt\rfloor-1)}(kx_1,x_2)^\top,
\eeqnn
where
$$
A_k=
\left( \begin{array}{ccc}
g'_{k,1,(1)}(1)&g'_{k,2,(1)}(1)\\
\\
g'_{k,1,(2)}(1)&g'_{k,2,(2)}(1)
\end{array}
\right ),\quad A_k^{(n)}=A_k\times A_k^{(n-1)},
$$
and $g'_{k,1,(i)}, g'_{k,2,(i)}, \Phi'_{k,1,(i)}$ and $\Phi'_{k,2,(i)}$ denote the derivative with respect to $\lambda_i, i=1,2$, respectively. By assumption $(\textbf{C,D})$, there exists $B\geq0$ such that $|\Phi'_{k,i,(j)}|\leq B$ for all $i,j\in\{1,2\}$ and $k\geq1$. It follows from Jensen's inequality that
\beqnn
\langle x,V_k(t,\theta)\rangle
\leq(\gamma_k^{-1}B+1)(\theta_1+\theta_2)(k^{-1},1)A_k^{(\lfloor\gamma_kt\rfloor-1)}(kx_1,x_2)^\top,
\eeqnn
then
\beqnn
&&\|V_k(t,\theta)\|=\sup_{\|x\|=1}\langle x,V_k(t,\theta)\rangle
\leq(\gamma_k^{-1}B+1)(\theta_1+\theta_2)\sup_{\|x\|=1}(k^{-1},1)B_k^{(\lfloor\gamma_kt\rfloor-1)}(kx_1,x_2)^\top\\
&&~~~~~~~~~~~~
\leq(2\gamma_k^{-1}B+1)^{\lfloor\gamma_kt\rfloor}\sqrt{2}(\theta_1+\theta_2),
\eeqnn
where
$$
B_k=
\left( \begin{array}{ccc}
\gamma_k^{-1}B+1&\gamma_k^{-1}Bk\\
\\
(k\gamma_k)^{-1}B&\gamma_k^{-1}B+1
\end{array}
\right ),
$$
by a modification of the proof of Lemma 2.6 and Theorem 2.7 in \cite{L19}, we get (\ref{3.8}). For given $(\Phi_1,\Phi_2)$ by (\ref{3.3})--(\ref{3.4}), it follows from Proposition \ref{pr:3.1} (ii) that there is a sequence $\{\Phi_{k,1},\Phi_{k,2}\}$ satisfying $(\textbf{A,C,D})$, let a sequence $\{V_k\}$ be given by (\ref{3.1}) and (\ref{3.2}), the existence of the solution is immediate. The uniqueness of the solution follows by Gronwall's inequality, and the semigroup property follows from the uniqueness of the solution.
\qed

\smallskip
\bproposition\label{pr:3.3}
Suppose $(\Phi_1,\Phi_2)$ are given by (\ref{3.3})--(\ref{3.4}), for any $\lambda\in\mbb{R}_{+}^2$ let $t\mapsto V(t,\lambda)$ be the unique positive solution to (\ref{3.8}). Then we can define a transition semigroup $(P_t)_{t\geq0}$ by
\beqlb
\int_\mbb{M}\mrm{e}^{-\langle\lambda,y\rangle}P_t(x,\mrm{d}y)=\exp\Big\{-\langle x,V(t,\lambda)\rangle\Big\},\quad x\in \mbb{M}.\label{3.10}
\eeqlb
\eproposition

\proof Given $(\Phi_1,\Phi_2)$ by (\ref{3.3})--(\ref{3.4}), by Proposition \ref{pr:3.1}, there is a sequence $(\Phi_{k,1},\Phi_{k,2})$ satisfying $(\textbf{A,C,D})$. By Proposition \ref{pr:3.2}, for any $a\geq0$ we have $V_k(t,\lambda)\rightarrow V(t,\lambda)$ uniformly on $[0,a]^3$ as $k\rightarrow\infty$. Taken $x_k\in\mbb{M}_k$ satisfying $x_k\rightarrow x$ as $k\rightarrow\infty,$
by a continuity theorem; see, e.g., Theorem 1.18 in \cite{L11}, (\ref{3.10}) defines a probability measure on $\mbb{M}$ and $\lim\limits_{k\rightarrow\infty}P_k^{\lfloor\gamma_kt\rfloor}(x_k,\cdot)=P_t(x,\cdot)$ by weak convergence. The semigroup property of the family of $(P_t)_{t\geq0}$ follows from (\ref{3.9}) and (\ref{3.10}).
\qed

\bdefinition
A Markov process $\{Y(t)=(Y_1(t),Y_2(t)):t\geq0\}$ is called a MSB-process with state space $\mbb{M},$ if it has the transition semigroup $(P_t)_{t\geq0}$ in (\ref{3.10}).
\edefinition

\bproposition\label{pr:3.4}
Let $(P_t)_{t\geq0}$ be the transition semigroup defined by (\ref{3.10}). Then we have
 \beqnn
\int_\mbb{M} \langle\lambda,y\rangle P_t(x,\mrm{d}y)= \langle x,\pi(t,\lambda)\rangle, \quad \lambda\in \mbb{R}_+^2, x\in \mbb{M},
 \eeqnn
where $t\mapsto \pi(t,\lambda)=(\pi_1(t,\lambda),\pi_2(t,\lambda))\in \mbb{R}_{+}^2$ is the unique solution to the following equations:
\beqnn
&&\frac{\mrm{d}\pi_1}{\mrm{d}t}(t,\lambda)= -a_{11}\pi_1(t,\lambda) + \pi_2(t,\lambda) \int_{\mbb{M}} z_2 n_1(\mrm{d}z), \cr
 \quad
&&\frac{\mrm{d}\pi_2}{\mrm{d}t}(t,\lambda)= a_{21}\pi_1(t,\lambda) + \int_{\mbb{M}_{-1}} \langle z,\pi(t,\lambda)\rangle n_2(\mrm{d}z)
\eeqnn
with initial condition $\pi(0,\lambda)=\lambda$.
\eproposition

\proof One can see that $V(t,0+)=0$ for $t\geq0$. By differentiating both sides of (\ref{3.10}) with respect to $\lambda_1$ and $\lambda_2,$ we have
\beqnn
&&\int_{\mbb{M}}y_1\,P_t(x,\mrm{d}y)=x_1\frac{\partial V_1(t,0+)}{\partial\lambda_1}+x_2\frac{\partial V_2(t,0+)}{\partial\lambda_1},\\ &&\int_{\mbb{M}}y_2\,P_t(x,\mrm{d}y)=x_1\frac{\partial V_1(t,0+)}{\partial\lambda_2}+x_2\frac{\partial V_2(t,0+)}{\partial\lambda_2}.
\eeqnn
It follows from (\ref{3.8}) that
\beqnn
&&\frac{\partial V_1(t,0+)}{\partial \lambda_1}=1+\int_0^t\Big(-a_{11}\frac{\partial V_1(s,0+)}{\partial \lambda_1}+\frac{\partial V_2(s,0+)}{\partial \lambda_1}\int_{\mbb{M}}z_2\,n_1(\mrm{d}z)\Big)\,\mrm{d}s,\\
&&\frac{\partial V_1(t,0+)}{\partial \lambda_2}=\int_0^t\Big(-a_{11}\frac{\partial V_1(s,0+)}{\partial \lambda_2}+\frac{\partial V_2(s,0+)}{\partial \lambda_2}\int_{\mbb{M}}z_2\,n_1(\mrm{d}z)\Big)\,\mrm{d}s,
\eeqnn
and ($\bar{a}_{21}:=a_{21}+\int_{\mbb{M}_{-1}}z_1\,n_2(\mrm{d}z)$)
\beqnn
&&\frac{\partial V_2(t,0+)}{\partial \lambda_1}=\int_0^t\Big(\frac{\partial V_1(s,0+)}{\partial \lambda_1}\bar{a}_{21}+\frac{\partial V_2(s,0+)}{\partial \lambda_1}\int_{\mbb{M}_{-1}}z_2\,n_2(\mrm{d}z)\Big)\,\mrm{d}s,\\
&&\frac{\partial V_2(t,0+)}{\partial \lambda_2}=1+\int_0^t\Big(\frac{\partial V_1(s,0+)}{\partial \lambda_2}\bar{a}_{21}+\frac{\partial V_2(s,0+)}{\partial \lambda_2}\int_{\mbb{M}_{-1}}z_2\,n_2(\mrm{d}z)\Big)\,\mrm{d}s.
\eeqnn
For $\theta=(\theta_1,\theta_2)\in\mbb{R}_{+}^2$ and $t\geq0$, we define $\pi(t,\theta)=(\pi_1(t,\theta),\pi_2(t,\theta))$ by
$$
\pi_1(t,\theta)=\theta_1\frac{\partial V_1(t,0+)}{\partial \lambda_1}+\theta_2\frac{\partial V_1(t,0+)}{\partial \lambda_2},\quad
\pi_2(t,\theta)=\theta_1\frac{\partial V_2(t,0+)}{\partial \lambda_1}+\theta_2\frac{\partial V_2(t,0+)}{\partial \lambda_2},
$$
we can conclude from the calculations above that
\beqnn
&&\pi_1(t,\theta)=\theta_1+\int_0^t\Big(-a_{11}\pi_1(s,\theta)+\int_{\mbb{M}}z_2\,n_1(\mrm{d}z)\pi_2(s,\theta)\Big)\,\mrm{d}s,\\
&&\pi_2(t,\theta)=\theta_2+\int_0^t\Big(a_{21}\pi_1(s,\theta)+\int_{\mbb{M}_{-1}}\langle z,\pi(s,\theta)\rangle\,n_2(\mrm{d}z)\Big)\,\mrm{d}s,
\eeqnn
and the desired assertion follows.
\qed

By a modification of the proof of Theorem 2.11 in \cite{L19}, one can see that the semigroup defined by (\ref{3.10}) is a Feller semigroup. Then the MSB-process has a $\mrm{c\grave{a}dl\grave{a}g}$ realization. Moreover, MSB-process can also be characterized in terms of a martingale problem described as follows, see Corollary 4.4 below for the proof. For $f\in C^2(\mbb{M}),$ let $L$ be an operator acting on $C^2(\mbb{M})$ defined by
\beqlb
&Lf(x)&= x_1\Big(\alpha f''_{11}(x)-a_{11}f'_1(x)+\int_\mbb{M}
\Big\{f(x+z)-f(x)-z_1f'_1(x)\Big\}\,n_1(\mrm{d}z)\Big)\nonumber\\
&&\quad\nonumber\\
&&~~~~~~~~+x_2\Big(a_{21}f'_1(x)+\int_{\mbb{M}_{-1}}\Big\{f(x+z)-f(x)\Big\}\,
n_2(\mrm{d}z)\Big).\nonumber
\eeqlb
Suppose that $\{(Y_1(t),Y_2(t)):t\geq0\}$ is a non-negative $\mrm{c\grave{a}dl\grave{a}g}$ process with $\mbf{E}[Y_i(0)]<\infty,i=1,2$. Then $\{(Y_1(t),Y_2(t)):t\geq0\}$ is a MSB-process with transition semigroup $(P_t)_{t\geq0}$ if and only if for every $f\in C^2(\mbb{M}),$
\beqnn
f(Y(t))=f(Y(0))+\int_0^tLf(Y(s))\,\mrm{d}s+\mrm{local}~\mrm{mart.}
\eeqnn

\btheorem\label{th:3.5}
Assume that $(\textbf{A,C,D})$ hold, $(Y_{k,1}(0)/k,Y_{k,2}(0))$ converges to $(Y_1(0),Y_2(0))$ in distribution. Then $\{(k^{-1}Y_{k,1}(\lfloor\gamma_kt\rfloor),Y_{k,2}(\lfloor\gamma_kt\rfloor)):t\geq0\}$ converges to $\{(Y_1(t),Y_2(t)):t\geq0\}$ in distribution on $D([0,\infty),\mbb{M})$ as $k\rightarrow\infty$.
\etheorem

\proof
Let $L$ be the generator of the MSB-process. For $\lambda=(\lambda_1,\lambda_2)\gg0$, $ x\in\mbb{M},$ set $\mrm{e}_{\lambda}(x)
=\mrm{e}^{-\langle\lambda,x\rangle}.$ We have
$$
L\mrm{e}_{\lambda}(x)=-\mrm{e}^{-\langle\lambda,x\rangle}\{x_1\Phi_1(\lambda)+x_2\Phi_2(\lambda)\}.
$$
Denote by $D_1$ the linear hull of $\{\mrm{e}_{\lambda},\lambda\gg0\}$, then $D_1$ is an algebra which strongly separates the points of $\mbb{M}$.
Let $C_0(\mbb{M})$ be the space of continuous function on $\mbb{M}$ vanishing at infinity. By the Stone-Weierstrass theorem, $D_1$ is dense in $C_0(\mbb{M})$ for the
supremum norm. Note that $D_1$ is invariant under $P_t$ by (\ref{3.10}), it follows from Proposition 3.3 in chapter I of \cite{EK86} that $D_1$ is the core of $L$. Note that $\{Y_{k,1}(n)/k,Y_{k,2}(n):n\geq0\}$ is a Markov chain with state space $\mbb{M}_k$, and one-step transition probability determined by:
$$
\int_{\mbb{M}_k}\mrm{e}^{-\langle\lambda,y\rangle}\,P_k(x,\mrm{d}y)=\Big(g_{k,1}(\mrm{e}^{-\lambda_1/k},\mrm{e}^{-\lambda_2})\Big)^{kx_1}\Big(g_{k,2}(\mrm{e}^{-\lambda_1/k},\mrm{e}^{-\lambda_2})\Big)^{x_2}.
$$
The (discrete) generator $L_k$ of $\{(k^{-1}Y_{k,1}(\lfloor\gamma_kt\rfloor),Y_{k,2}(\lfloor\gamma_kt\rfloor)):t\geq0\}$ is given by:
\beqnn
&&L_k\mrm{e}_{\lambda}(x)=\gamma_k\Big\{\Big(g_{k,1}(\mrm{e}^{-\lambda_1/k},\mrm{e}^{-\lambda_2})\Big)^{kx_1}\cdot\Big(g_{k,2}
(\mrm{e}^{-\lambda_1/k},\mrm{e}^{-\lambda_2})\Big)^{x_2}-\mrm{e}^{-\langle\lambda,x\rangle}\Big\}\\
&&~~~~~~
=\mrm{e}^{-\langle\lambda,x\rangle}\gamma_k\Big\{\exp[kx_1\log(1-(k\gamma_k)^{-1}\Phi_{k,1}(\lambda_1,\lambda_2)\mrm{e}^{\lambda_1/k})\\
&&~~~~~~~~~~~~~~~~~~~~~~~~~~~~
+x_2\log(1-\gamma_k^{-1}\Phi_{k,2}(\lambda_1,\lambda_2)\mrm{e}^{\lambda_2})]-1\Big\}\\
&&~~~~~~
=\mrm{e}^{-\langle\lambda,x\rangle}\Big\{x_1\bar{\Phi}_{k,1}(\lambda_1,\lambda_2)+x_2\bar{\Phi}_{k,2}(\lambda_1,\lambda_2)\Big\}+o(1),
\eeqnn
where $\bar{\Phi}_{k,i},\Phi_{k,i}, i=1,2$ are defined as before. It follows from Proposition \ref{pr:3.1} that
$$\lim\limits_{k\rightarrow\infty}\sup\limits_{x\in E_k}|L_k\mrm{e}_{\lambda}(x)-L\mrm{e}_{\lambda}(x)|=0.
$$
From Corollary 8.9 in Chapter 4 of \cite{EK86}, we prove the desired result.
\qed

\btheorem\label{th:3.6}
Suppose that $\{(Y_1(t),Y_2(t)):t\geq0\}$ is any MSB-process with branching mechanism $(\Phi_1,\Phi_2)$. Then there exist a sequence of positive numbers $\{\gamma_k\}$ and a sequence of two-type GW-processes
$\{(Y_{k,1}(n),Y_{k,2}(n)):n\in\mbb{N}\}$ with generation functions $(g_{k,1},g_{k,2})$ such that the sequence
$\{(k^{-1}Y_{k,1}(\lfloor\gamma_kt\rfloor),$
$Y_{k,2}(\lfloor\gamma_kt\rfloor)):t\geq0\}$ converges in distribution on $D([0,\infty),\mbb{M})$ to the process $\{(Y_1(t),Y_2(t)):t\geq0\}$ as $k\rightarrow\infty$.
\etheorem

\proof
By Proposition \ref{pr:3.1}, there exist $\{\gamma_k\},\{(g_{k,1},g_{k,2})\}$ such that conditions $(\textbf{A,C,D})$ hold. And the desired result follows from Theorem \ref{th:3.5}.
\qed

\section{The construction of MSB-processes by stochastic equations}
\setcounter{equation}{0}
Let $(\Omega,{\cal F},{\cal F}_{t},\mbf{P})$ be a complete filtered probability space satisfying the usual hypotheses, let $\{B(t)\}$ be a standard Brownian motion, $\{N_1(\mrm{d}s,\mrm{d}u,\mrm{d}z)\}$ be a Poisson random measure on $(0,\infty)^2\times\mbb{M}$ with intensity
$ \mrm{d}s\mrm{d}un_1(\mrm{d}z)$, $\{N_2(\mrm{d}s,\mrm{d}u,\mrm{d}z)\}$ be a Poisson random measure on $(0,\infty)^2\times {\mbb{M}_{-1}}$ with intensity $\mrm{d}s\mrm{d}un_2(\mrm{d}z)$, $ z=(z_1,z_2).$ Suppose that $B,N_1,N_2$ are independent of each other. Let us recall the  stochastic integral equation system (\ref{1.7})--(\ref{1.8}):
\beqnn
&&Y_1(t)=Y_1(0)-\int_0^ta_{11}Y_1(s)\,\mrm{d}s+\int_0^t\sqrt{2\alpha Y_1(s)}\,\mrm{d}B(s)+\int_0^t\int_0^{Y_1(s-)}\int_\mbb{M}z_1\,
\tilde{N}_1(\mrm{d}s,\mrm{d}u,\mrm{d}z)\\
&&\qquad~~~~~~~~~~~
+\int_0^ta_{21}Y_2(s)\,\mrm{d}s+\int_0^t\int_0^{Y_2(s-)}\int_{\mbb{M}_{-1}}z_1\,N_2(\mrm{d}s,\mrm{d}u,\mrm{d}z), \\
&&\qquad\\
&&Y_2(t)=Y_2(0)+\int_0^t\int_0^{Y_1(s-)}\int_\mbb{M}z_2\,N_1(\mrm{d}s,\mrm{d}u,\mrm{d}z)
+\int_0^t\int_0^{Y_2(s-)}\int_{\mbb{M}_{-1}}z_2\,N_2(\mrm{d}s,\mrm{d}u,\mrm{d}z),
\eeqnn
where $\tilde{N}_1(\mrm{d}s,\mrm{d}u,\mrm{d}z)=N_1(\mrm{d}s,\mrm{d}u,\mrm{d}z)-\mrm{d}s\mrm{d}un_1(\mrm{d}z)$ is the compensated Poisson random measure of $N_1$.

\bproposition\label{pr:4.1}
Suppose that $\{Y(t)\}$ satisfies (\ref{1.7})--(\ref{1.8}) and $\mbf{P}\{Y(0)\geq\mbf{0}\}=1$. Then $\mbf{P}\{Y(t)\geq\mbf{0},\forall t\geq0\}=1.$
\eproposition

 \proof
By equation (\ref{1.8}), if $Y_2(0)\geq0,$ it is not hard to see that for all $t\geq0,$ $Y_2(t)\geq0.$ Now suppose that there exists $\varepsilon>0,$ such that $\tau:=\inf\{t>0,Y_1(t)\leq-\varepsilon\}<\infty$ with strictly positive probability. Then there exists $t_0>0,Y_1(t_0)=0$, and on the time interval $[t_0,\tau], t\mapsto Y_1(t)$ is a strictly negative continuous function. Hence there are some $t_1\in[t_0,\tau]$ and $\delta>0,$ such that for all $s\in[t_0,t_1], -a_{11}Y_1(s)+a_{21}Y_2(s)\geq\delta.$ Then
$$
Y_1(t_1)=Y_1(t_1)-Y_1(t_0)\geq\int_{t_0}^{t_1}(-a_{11}Y_1(s)+a_{21}Y_2(s))\,\mrm{d}s\geq\delta(t_1-t_0)>0,
$$
since $Y_1(t)<0,\forall t\in(t_0,\tau]$, we get a contradiction.
\qed

To analyze the property of above equation system, we first construct a sequence of functions $\{\phi_k\}$ on $\mbb{R}$ as follows. For each integer
$k\geq0$ define $a_k=\exp\{-k(k+1)/2\}$. Then $a_k\rightarrow0$ decreasingly as $k\rightarrow\infty$ and $\int_{a_k}^{a_{k-1}}z^{-1}\,\mrm{d}z=k$ for $k\geq1$. Let $x\mapsto\psi_k(x)$ be a non-negative continuous function on $\mbb{R}$ which has support in $(a_k,a_{k-1})$ and satisfies $\int_{a_k}^{a_{k-1}}\psi_k(x)\,\mrm{d}x=1$ and $0\leq\psi_k(x)\leq2(kx)^{-1}$ for $a_k<x<a_{k-1}.$ For each $k\geq1$ let
\beqlb
\phi_k(z)=\int_0^{|z|}\,\mrm{d}y\int_0^{y}\psi_k(x)\,\mrm{d}x,\quad z\in\mbb{R}.\nonumber
\eeqlb
Moreover,  for function $f$ on $\mbb{R}$ we denote
$$
\Delta_{z}f(x)=f(x+z)-f(x).
$$

\btheorem\label{th:4.2}
The pathwise uniqueness for (\ref{1.7})--(\ref{1.8}) holds.
\etheorem

 \proof
Suppose that $\{Y(t)\}$ and $\{Y'(t)\}$ are two solutions of (\ref{1.7})--(\ref{1.8}). Let $\zeta_i(t)=Y_i(t)-Y'_i(t), i=1,2$ for $t\geq0.$ We have
\beqnn
&&\zeta_1(t)=\zeta_1(0)-\int_0^t\Big(a_{11}\zeta_1(s)+a_{21}\zeta_2(s)\Big)\,\mrm{d}s
+\int_0^t\Big(\sqrt{2\alpha Y_1(s)}-\sqrt{2\alpha Y'_1(s)}\Big)\,\mrm{d}B(s)\\
&&\qquad~~~~~~~~~~
+\int_0^t\int_{Y'_1(s-)}^{Y_1(s-)}\int_\mbb{M} z_1\mbf{1}_{\{\zeta_1(s-)>0\}}\,\tilde{N}_1(\mrm{d}s,\mrm{d}u,\mrm{d}z)\\
&&\qquad~~~~~~~~~~
-\int_0^t\int_{Y_1(s-)}^{Y'_1(s-)}\int_\mbb{M} z_1\mbf{1}_{\{\zeta_1(s-)\leq0\}}\,\tilde{N}_1(\mrm{d}s,\mrm{d}u,\mrm{d}z)\\
&&\qquad~~~~~~~~~~
+\int_0^t\int_{Y'_2(s-)}^{Y_2(s-)}\int_{\mbb{M}_{-1}} z_1\mbf{1}_{\{\zeta_2(s-)>0\}}\,N_2(\mrm{d}s,\mrm{d}u,\mrm{d}z)\\
&&\qquad~~~~~~~~~~
-\int_0^t\int_{Y_2(s-)}^{Y'_2(s-)}\int_{\mbb{M}_{-1}} z_1\mbf{1}_{\{\zeta_2(s-)\leq0\}}\,N_2(\mrm{d}s,\mrm{d}u,\mrm{d}z).
\eeqnn
Let $\tau_m=\inf\Big\{t\geq0:Y_1(t)\vee Y_2(t) \vee Y'_1(t)\vee Y'_2(t)\geq m\Big\}$ for $m\geq1.$ By similar calculations as in Theorem 3.4 in \cite{M14}, there exists $C_1>0$ such that
$$
\mbf{E}|\zeta_1(t\wedge\tau_m)|\leq C_1\int_0^t\mbf{E}\Big(|\zeta_1(s\wedge\tau_m)|
+|\zeta_2(s\wedge\tau_m)|\Big)\,\mrm{d}s.
$$
On the other hand,
\beqnn
&&\zeta_2(t)=\zeta_2(0)+\int_0^t\int_{Y'_1(s-)}^{Y_1(s-)}\int_\mbb{M} z_2\mbf{1}_{\{\zeta_1(s-)>0\}}\,N_1(\mrm{d}s,\mrm{d}u,\mrm{d}z)\\
&&\qquad~~~~~~~~~~
-\int_0^t\int_{Y_1(s-)}^{Y'_1(s-)}\int_\mbb{M} z_2\mbf{1}_{\{\zeta_1(s-)\leq0\}}\,N_1(\mrm{d}s,\mrm{d}u,\mrm{d}z)\\
&&\qquad~~~~~~~~~~
+\int_0^t\int_{Y'_2(s-)}^{Y_2(s-)}\int_{\mbb{M}_{-1}} z_2\mbf{1}_{\{\zeta_2(s-)>0\}}\,N_2(\mrm{d}s,\mrm{d}u,\mrm{d}z)\\
&&\qquad~~~~~~~~~~
-\int_0^t\int_{Y_2(s-)}^{Y'_2(s-)}\int_{\mbb{M}_{-1}} z_2\mbf{1}_{\{\zeta_2(s-)\leq0\}}\,N_2(\mrm{d}s,\mrm{d}u,\mrm{d}z).
\eeqnn
By $\mrm{It\hat{o}}$'s formula,
\beqnn
&&\phi_k(\zeta_2(t\wedge\tau_m))=\phi_k(\zeta_2(0))+\int_0^{t\wedge\tau_m}\int_{Y'_1(s-)}^{Y_1(s-)}\int_\mbb{M}\Delta_{z_2}\phi_k(\zeta_2(s-))\mbf{1}_{\{\zeta_1(s-)>0\}}\,
\mrm{d}s\,\mrm{d}u\,n_1(\mrm{d}z)\\
&&\qquad~~~~~~~~~~~~~~~~
+\int_0^{t\wedge\tau_m}\int_{Y_1(s-)}^{Y'_1(s-)}\int_\mbb{M}\Delta_{-z_2}\phi_k(\zeta_2(s-))\mbf{1}_{\{\zeta_1(s-)\leq0\}}\,
\mrm{d}s\,\mrm{d}u\,n_1(\mrm{d}z)\\
&&\qquad~~~~~~~~~~~~~~~~
+\int_0^{t\wedge\tau_m}\int_{Y'_2(s-)}^{Y_2(s-)}\int_{\mbb{M}_{-1}}\Delta_{z_2}\phi_k(\zeta_2(s-))\mbf{1}_{\{\zeta_2(s-)>0\}}\,
\mrm{d}s\,\mrm{d}u\,n_2(\mrm{d}z)\\
&&\qquad~~~~~~~~~~~~~~~~
+\int_0^{t\wedge\tau_m}\int_{Y_2(s-)}^{Y'_2(s-)}\int_{\mbb{M}_{-1}}\Delta_{-z_2}\phi_k(\zeta_2(s-))\mbf{1}_{\{\zeta_2(s-)\leq0\}}\,
\mrm{d}s\,\mrm{d}u\,n_2(\mrm{d}z)+\mrm{mart}.
\eeqnn
Similarly, there exists $C_2>0$ such that
$$
\mbf{E}|\zeta_2(t\wedge\tau_m)|\leq C_2\int_0^t\mbf{E}\Big(|\zeta_1(s\wedge\tau_m)|
+|\zeta_2(s\wedge\tau_m)|\Big)\,\mrm{d}s.
$$
In conclusion,
$$
\mbf{E}\Big(|\zeta_1(t\wedge\tau_m)|+|\zeta_2(t\wedge\tau_m)|\Big)\leq (C_1+C_2)\int_0^t\mbf{E}\Big(|\zeta_1(s\wedge\tau_m)|+|\zeta_2(s\wedge\tau_m)|\Big)\,\mrm{d}s.
$$
By Gronwall inequality, for all $t\geq0,$
$$
\mbf{E}\Big(|\zeta_1(t\wedge\tau_m)|+|\zeta_2(t\wedge\tau_m)|\Big)=0.
$$
Since $\{Y(t)\}$ and $\{Y'(t)\}$ have $\mrm{c\grave{a}dl\grave{a}g}$ sample paths, we conclude that $\mbf{P}\{Y(t)=Y'(t),\forall t\geq0\}=1$ as $m\rightarrow\infty.$
\qed

\btheorem\label{co:4.3}
There is a unique non-negative strong solution to (\ref{1.7})--(\ref{1.8}).
\etheorem

 \proof
Since $\nu_1$ is supported on $\mbb{M}\setminus\{\mbf{0}\},$ we can rewrite (\ref{1.7})--(\ref{1.8}) as
\beqlb
&&Y_1(t)=Y_1(0)+\int_0^t\Big(a_{21}Y_2(s)-a_{11}Y_1(s)\Big)\,\mrm{d}s+\int_0^t\sqrt{2\alpha Y_1(s)}\,\mrm{d}B(s)
\nonumber \\
&&\qquad~~~~~~~~~~~
+\int_0^t\int_0^{Y_1(s-)}\int_{\mbb{M}_{-1}}z_1\,\tilde{N}_1(\mrm{d}s,\mrm{d}u,\mrm{d}z)+\int_0^t\int_0^{Y_2(s-)}\int_{\mbb{M}_{-1}}z_1\,N_2(\mrm{d}s,\mrm{d}u,\mrm{d}z), \nonumber\\
&&\qquad \nonumber\\
&&Y_2(t)=Y_2(0)+\int_0^t\int_0^{Y_1(s-)}\int_{\mbb{M}_{-1}}z_2\,N_1(\mrm{d}s,\mrm{d}u,\mrm{d}z)+\int_0^t\int_0^{Y_2(s-)}\int_{\mbb{M}_{-1}}z_2\,N_2(\mrm{d}s,\mrm{d}u,\mrm{d}z).\nonumber
\eeqlb
For any fixed $n\geq1,$ let $V_n=\{z\in {\mbb{M}_{-1}}:\|z\|\geq1/n\},$ then $n_1(V_n)+n_2(V_n)<\infty.$ For $m\geq1$ and $x\in\mbb{M}$, define
\beqnn
&&b(x,m)=a_{21}(x_2\wedge m)-a_{11}(x_1\wedge m),\quad \theta(m,n)=\int_{V_n}(z_1\wedge m)\,n_1(\mrm{d}z),\\
&&\beta_1(m)=\int_{\mbb{M}_{-1}}(z_1-z_1\wedge m)\,n_1(\mrm{d}z),\quad \beta_2(m)=\int_{\mbb{M}_{-1}}(z_2-z_2\wedge m)\,n_2(\mrm{d}z).
\eeqnn
By the results for continuous-type stochastic equations in \cite[p.169]{IW89}, one can show that there is a non-negative weak solution to the following
stochastic equation system:
\beqnn
&&Y_1(t)=Y_1(0)+\int_0^t\Big(b(Y(s),m)-(\beta_1(m)+\theta(m,n))(Y_1(s)\wedge m)\Big)\,\mrm{d}s\\
&&~~~~~~~~~~~~~~~~~
+\int_0^t\sqrt{2\alpha(Y_1(s)\wedge m)}\,\mrm{d}B(s),\\
&&Y_2(t)=Y_2(0)-\int_0^t\beta_2(m)(Y_2(s)\wedge m)\,\mrm{d}s.
\eeqnn
The pathwise uniqueness holds for the above system of equations by similar arguments as in Theorem \ref{th:4.2}. Then it has a unique strong solution. By similar arguments as in the proof of Proposition 2.2 in \cite{FL10}, we can get a pathwise unique non-negative strong solution $\{Y_{m,n}(t):t\geq0\}$ to (\ref{4.1})--(\ref{4.2}):
\beqlb
&&Y_1(t)=Y_1(0)+\int_0^tb(Y(s),m)-\beta_1(m)(Y_1(s)\wedge m)\,\mrm{d}s+\int_0^t\sqrt{2\alpha (Y_1(s)\wedge m)}\,\mrm{d}B(s)\nonumber \\
&&\qquad~~~~~~~~~~~
+\int_0^t\int_0^{Y_1(s-)\wedge m}\int_{V_n}(z_1\wedge m)\,\tilde{N}_1(\mrm{d}s,\mrm{d}u,\mrm{d}z)\nonumber\\
&&\qquad~~~~~~~~~~~+\int_0^t\int_0^{Y_2(s-)\wedge m}\int_{V_n}(z_1\wedge m)\,N_2(\mrm{d}s,\mrm{d}u,\mrm{d}z),\label{4.1}\\
&&\qquad \nonumber\\
&&Y_2(t)=Y_2(0)-\int_0^t\beta_2(m)(Y_2(s)\wedge m)\,\mrm{d}s+\int_0^t\int_0^{Y_1(s-)\wedge m}\int_{V_n}(z_2\wedge m)\,N_1(\mrm{d}s,\mrm{d}u,\mrm{d}z)\nonumber\\
&&\qquad~~~~~~~~~~~+\int_0^t\int_0^{Y_2(s-)\wedge m}\int_{V_n}(z_2\wedge m)\,N_2(\mrm{d}s,\mrm{d}u,\mrm{d}z).\label{4.2}
\eeqlb
As in the proof of Lemma 4.3 in \cite{FL10}, one can see the sequence $ \{Y_{m,n}(t):t\geq0\},n=1,2,...$ is tight in $D([0,\infty),\mbb{M}).$ Following the proof of Theorem 4.4 in \cite{FL10}, it is easy to show that any weak limit point $\{Y_m(t):t\geq0\}$ of the sequence is a non-negative weak solution to
\beqlb
&&Y_1(t)=Y_1(0)+\int_0^tb(Y(s),m)-\beta_1(m)(Y_1(s)\wedge m)\,\mrm{d}s+\int_0^t\sqrt{2\alpha (Y_1(s)\wedge m)}\,\mrm{d}B(s)\nonumber\\
&&\qquad~~~~~~~~~~~
+\int_0^t\int_0^{Y_1(s-)\wedge m}\int_{\mbb{M}_{-1}}(z_1\wedge m)\,\tilde{N}_1(\mrm{d}s,\mrm{d}u,\mrm{d}z)\nonumber\\
&&\qquad~~~~~~~~~~~
+\int_0^t\int_0^{Y_2(s-)\wedge m}\int_{\mbb{M}_{-1}}(z_1\wedge m)\,N_2(\mrm{d}s,\mrm{d}u,\mrm{d}z),\label{4.3} \\
&&\qquad \nonumber\\
&&Y_2(t)=Y_2(0)-\int_0^t\beta_2(m)(Y_2(s)\wedge m)\,\mrm{d}s+\int_0^t\int_0^{Y_1(s-)\wedge m}\int_{\mbb{M}_{-1}}(z_2\wedge m)\,
N_1(\mrm{d}s,\mrm{d}u,\mrm{d}z)\nonumber\\
&&\qquad~~~~~~~~~~~
+\int_0^t\int_0^{Y_2(s-)\wedge m}\int_{\mbb{M}_{-1}}(z_2\wedge m)\,N_2(\mrm{d}s,\mrm{d}u,\mrm{d}z).\label{4.4}
\eeqlb
By Theorem \ref{th:4.2}, the pathwise uniqueness holds for (\ref{4.3})--(\ref{4.4}), so the system of equations has a unique strong solution. Finally, the desired result
follows from a modification of the proof of Proposition 2.4 in \cite{FL10}.
\qed

\smallskip

\bcorollary\label{co:4.4}
A $c\grave{a}dl\grave{a}g$ non-negative process is a MSB-process with transition semigroup $(P_t)_{t\geq0}$ defined by (\ref{3.8}) and (\ref{3.10}) if and only if it is a weak solution of (\ref{1.7})--(\ref{1.8}).
\ecorollary

 \proof
 Suppose that $\{(Y_1(t),Y_2(t))\}_{t\geq0}$ is a weak solution of (\ref{1.7})--(\ref{1.8}). By $\mathrm{It\hat{o}}$'s formula one can see that $\{(Y_1(t),Y_2(t))\}_{t\geq0}$ solves the martingale problem associated with the generator $L$. By the arguments in Section 3 we infer that $\{(Y_1(t),Y_2(t))\}_{t\geq0}$ is a MSB-process with transition semigroup $(P_t)_{t\geq0}$ defined by (\ref{3.8}) and (\ref{3.10}). Conversely, suppose that $\{(Y_1(t),Y_2(t))\}_{t\geq0}$ is a $\mrm{c\grave{a}dl\grave{a}g}$ realization of MSB-process with transition semigroup $(P_t)_{t\geq0}$ defined by (\ref{3.8}) and (\ref{3.10}). Then the distributions of $\{(Y_1(t),Y_2(t))\}_{t\geq0}$ on $D([0,\infty),\mbb{M})$ can be characterized uniquely by the martingale problem. By a standard stopping time argument, we have
\beqlb
&&Y_1(t)=Y_1(0)-\int_0^t(a_{11}Y_1(s)-a_{21}Y_2(s))\,\mrm{d}s
+\int_0^t\int_{\mbb{M}_{-1}}Y_2(s)z_1\,\mrm{d}s\,n_2(\mrm{d}z)+G_1(t),\label{4.5}\\
&&Y_2(t)=Y_2(0)+\int_0^t\int_{\mbb{M}}Y_1(s)z_2\,\mrm{d}s\,n_1(\mrm{d}z)
+\int_0^t\int_{\mbb{M}_{-1}}Y_2(s)z_2\,\mrm{d}s\,n_2(\mrm{d}z)+G_2(t),\label{4.6}
\eeqlb
where $G_1(t)$ and $G_2(t)$ are two square-integrable local martingales. Let $N_0(\mrm{d}s,\mrm{d}z)$ be the optimal random measure on $[0,\infty)\times {\mbb{M}_{-1}}$ defined by
\beqlb
N_0(\mrm{d}s,\mrm{d}z):=\sum\limits_{s>0}\mbf{1}_{\{(Y_1(s),Y_2(s))\neq(Y_1(s-),Y_2(s-))\}}\delta_{(s,Y(s)-Y(s-))}(\mrm{d}s,\mrm{d}z),\nonumber
\eeqlb
It follows from \cite[p.376]{DM82} that
\beqnn
G_1(t)=G^c_1(t)+\int_0^t\int_{\mbb{M}_{-1}}z_1\,\tilde{N}_0(\mrm{d}s,\mrm{d}z),\quad
G_2(t)=G^c_2(t)+\int_0^t\int_{\mbb{M}_{-1}}z_2\,\tilde{N}_0(\mrm{d}s,\mrm{d}z),
\eeqnn
where $t\mapsto G^c_1(t)$ and $t\mapsto G^c_2(t)$ are two continuous local martingales with quadratic variations $t\mapsto C_1(t)$ and $t\mapsto C_2(t),$ respectively. Using $\mathrm{It\hat{o}}$'s formula to (\ref{4.5})--(\ref{4.6}) and the uniqueness of canonical decompositions of semi-martingales we find that
$N_0(\mrm{d}s,\mrm{d}z)$ has a predictable compensator
\beqlb
\hat{N}_0(\mrm{d}s,\mrm{d}z)=Y_1(s-)\,\mrm{d}s\,n_1(\mrm{d}z)+Y_2(s-)\,\mrm{d}s\,n_2(\mrm{d}z),\nonumber
\eeqlb
$\mrm{d}C_1(t)=2\alpha Y_1(t)\,\mrm{d}t$ and $\mrm{d}C_2(t)=0$. Then we obtain the equation
(\ref{1.7})--(\ref{1.8}) on an extension of the probability space by applying martingale representation theorems; see, e.g., \cite[p.93, p.84]{IW89}, which completes the proof.
\qed

\section{The distribution of local jumps}
\setcounter{equation}{0}For any initial time $r\geq0$, let $Y=(\Omega,{\cal F},{\cal F}_{r,t},Y(t), \mbf{P}_{r,y}:t\geq r,y\geq0)$ be a Hunt realization of the MSB-process with transition semigroup $(P_t)_{t\geq0}$ defined by $(\ref{3.8})$ and $(\ref{3.10})$. Here, $\{\mbf{P}_{r,y}:y\geq0\}$ be a family of probability measures on $(\Omega,{\cal F},{\cal F}_{r,t})$ satisfying $\mbf{P}_{r,y}\{Y(r)=y\}=1$ for all $y\geq0$. For any $t\geq r\geq0$ and $\lambda\in[0,\infty)^2,$ we have
\beqnn
\mbf{P}_{r,y}\,\mrm{exp}\Big\{-\langle\lambda,Y(t)\rangle\Big\}=\mrm{exp}\Big\{-\langle
y,\bar{V}(r,\lambda)\rangle\Big\},
\eeqnn
where $r\rightarrow \bar{V}(r,\lambda):= V(t-r,\lambda)$ satisfies
\beqnn
\bar{V}_1(r,\lambda)=\int_r^t \Phi_1(\bar{V}(s,\lambda))\,\mrm{d}s+\lambda_1,\quad \bar{V}_2(r,\lambda)=\int_r^t \Phi_2(\bar{V}(s,\lambda))\,\mrm{d}s+\lambda_2,\quad 0\leq r\leq t.
\eeqnn
By modifying the arguments of Proposition 4.1, Theorem 4.2 and Corollary 4.4 in \cite{L19}, we have the following

\bproposition\label{pr:5.1}
For $\{t_1<...<t_n\}\subset[0,\infty)$ and $\{\lambda_1,...\lambda_n\}\subset[0,\infty)^2$, we have
\beqnn
\mbf{P}_{r,y}\,\mrm{exp}\Big\{-\sum_{j=1}^n\langle\lambda_j,Y(t_j)\mbf{1}_{\{r\leq t_j\}}\rangle\Big\}
=\mrm{exp}\Big\{-\langle y,\bar{V}(r)\rangle\Big\},\quad 0\leq r\leq t_n,
\eeqnn
where $\bar{V}(r)=\bar{V}(r,\lambda_1,...,\lambda_n)$ on $[0,t_n]$ satisfies
$$
\bar{V}_1(r)=\int_r^{t_n}\Phi_1(\bar{V}(s))\,\mrm{d}s+\sum_{j=1}^n\lambda_{j1}\mbf{1}_{\{r\leq t_j\}},\quad \bar{V}_2(r)=\int_r^{t_n}\Phi_2(\bar{V}(s))\,\mrm{d}s+\sum_{j=1}^n\lambda_{j2}\mbf{1}_{\{r\leq t_j\}}.
$$
\eproposition

\bproposition\label{pr:5.2}
Suppose that  $t\geq0$ and $\mu$ is a finite measure supported by $[0,t].$ Let $s\mapsto\lambda(s)=(\lambda_1(s),\lambda_2(s))$ be a bounded positive Borel function on $[0,t],$ then we have
$$
\mbf{P}_{r,y}\mrm{exp}\Big\{-\int_{[r,t]}\langle\lambda(s),Y(s)\rangle\,\mu(\mrm{d}s)\Big\}
=\mrm{exp}\Big\{-\langle y,\bar{V}(r)\rangle\Big\},\quad 0\leq r\leq t,
$$
where $r\mapsto\bar{V}(r)=\bar{V}(r,\lambda_1,\lambda_2)$ is the positive solution on $[0,t]$ of
$$
\bar{V}_1(r)=\int_r^{t}\Phi_1(\bar{V}(s))\,\mrm{d}s+\int_r^t\lambda_1(s)\,\mu(\mrm{d}s),\quad \bar{V}_2(r)=\int_r^{t}\Phi_2(\bar{V}(s))\,\mrm{d}s+\int_r^t\lambda_2(s)\,\mu(\mrm{d}s).
$$
\eproposition

\bcorollary\label{co:5.3}
Let $Y=(\Omega,{\cal F},{\cal F}_t,Y(t),\mbf{P}_y)$ be a Hunt realization of the MSB-process started from time zero. Then for $t\geq0, \lambda=(\lambda_1,\lambda_2)\in[0,\infty)^2,$ we have
\beqlb\label{5.1}
\mbf{P}_y\,\mrm{exp}\Big\{-\int_0^t\langle\lambda,Y(s)\rangle\,\mrm{d}s\Big\}=\mrm{exp}\Big\{-\langle y,\bar{V}(t)\rangle\Big\},
\eeqlb
where $t\mapsto\bar{V}(t)=\bar{V}(t,\lambda)$ is the positive solution of
\beqlb
\bar{V}_1(t)=\int_0^{t}\Phi_1(\bar{V}(s))\,\mrm{d}s+\lambda_1t,\quad \bar{V}_2(t)=\int_0^{t}\Phi_2(\bar{V}(s))\,\mrm{d}s+\lambda_2t.
\eeqlb
\ecorollary

We shall introduce some notations before presenting the main results in this section. Let $r=(r_1,r_2)\in[0,\infty)^2,~ A_r=(r_1,\infty)\times(r_2,\infty)$ and
$n(A_r)=(n_1(A_r),n_2(A_r)).$ We define two functions on $[0,\infty)^2$:
\beqlb
&&\Phi_1^r(\lambda_1,\lambda_2)=-a_{11}^{r}\lambda_1+b_{11}^{r}\lambda_2-\alpha\lambda_1^2-\int_{\mbb{M}\setminus A_r}(\mrm{e}^{-\langle\lambda,z\rangle}-1+\langle\lambda,z\rangle)\,n_1(\mrm{d}z),\label{5.3}\\
&&\Phi_2^r(\lambda_1,\lambda_2)=a_{21}^{r}\lambda_1+b_{21}^{r}\lambda_2-\int_{{\mbb{M}_{-1}}\setminus A_r}(\mrm{e}^{-\langle\lambda,z\rangle}-1+\langle\lambda,z\rangle)\,n_2(\mrm{d}z),\label{5.4}
\eeqlb
where
\beqlb
&&a_{11}^{r}=a_{11}+\int_{A_r}z_1\,n_1(\mrm{d}z),\quad a_{21}^{r}=a_{21}+\int_{{\mbb{M}_{-1}}\setminus A_r}z_1\,n_2(\mrm{d}z), \nonumber\\
&&\quad \nonumber\\
&&b_{11}^{r}=\int_{\mbb{M}\setminus A_r}z_2\,n_1(\mrm{d}z),\quad b_{21}^{r}=\int_{{\mbb{M}_{-1}}\setminus A_r}z_2\,n_2(\mrm{d}z).\nonumber
\eeqlb

The following theorem gives a characterization of the distribution of the local maximal jump of the MSB-process.

\smallskip

\btheorem\label{th:5.4}
  Let $\tau_r=\inf\{s\geq0:\Delta Y_1(s)>r_1~\mrm{or}~\Delta Y_2(s)>r_2 \}$. Then we have
\beqlb
\mbf{P}_y(\tau_r>t)=\mrm{exp}\Big\{-\langle y,\bar{V}^r(t,n(A_r))\rangle\Big\},\nonumber
\eeqlb
where $\bar{V}^r(t)=(\bar{V}_1^r(t),\bar{V}_2^r(t))$ is the solution of
$$
\bar{V}_1^r(t)=\int_0^{t}\Phi_1^r(\bar{V}^r(s))\,\mrm{d}s+n_1(A_r)t,\quad \bar{V}_2^r(t)=\int_0^{t}\Phi_2^r(\bar{V}^r(s))\,\mrm{d}s+n_2(A_r)t.
$$
\etheorem

\smallskip
 \proof
We can rewrite equations (\ref{1.7})--(\ref{1.8}) by:
\beqnn
&Y_1(t)&=Y_1(0)+\int_0^t\Big(a_{21}^{r}Y_2(s)-a_{11}^{r}Y_1(s)\Big)\,\mrm{d}s+\int_0^t\sqrt{2\alpha Y_1(s)}\,\mrm{d}B(s)\\
&&\qquad\\
&&+\int_0^t\int_0^{Y_1(s-)}\int_{\mbb{M}\setminus A_r}z_1\,\tilde{N}_1(\mrm{d}s,\mrm{d}u,\mrm{d}z)+\int_0^t\int_0^{Y_2(s-)}\int_{{\mbb{M}_{-1}}\setminus A_r}z_1\,\tilde{N}_2(\mrm{d}s,\mrm{d}u,\mrm{d}z)\\
&&\qquad\\
&&+\int_0^t\int_0^{Y_1(s-)}\int_{A_r}z_1\,N_1(\mrm{d}s,\mrm{d}u,\mrm{d}z)+\int_0^t\int_0^{Y_2(s-)}\int_{A_r}z_1\,N_2(\mrm{d}s,\mrm{d}u,\mrm{d}z),\\
&&\qquad\\
&Y_2(t)&=Y_2(0)+\int_0^t\Big(b_{11}^{r}Y_1(s)+b_{21}^{r}Y_2(s)\Big)\,\mrm{d}s\\
&&\qquad\\
&&+\int_0^t\int_0^{Y_1(s-)}\int_{\mbb{M}\setminus A_r}z_2\,\tilde{N}_1(\mrm{d}s,\mrm{d}u,\mrm{d}z)+\int_0^t\int_0^{Y_2(s-)}\int_{{\mbb{M}_{-1}}\setminus A_r}z_2\,\tilde{N}_2(\mrm{d}s,\mrm{d}u,\mrm{d}z)\\
&&\qquad\\
&&+\int_0^t\int_0^{Y_1(s-)}\int_{A_r}z_2\,N_1(\mrm{d}s,\mrm{d}u,\mrm{d}z)+\int_0^t\int_0^{Y_2(s-)}\int_{A_r}z_2\,N_2(\mrm{d}s,\mrm{d}u,\mrm{d}z).
\eeqnn
Let $\Big((Y_1^{r_1}(t),Y_2^{r_2}(t)):t\geq0\Big)$ be the strong solution to
\beqnn
&Y_1^{r_1}(t)&=Y_1(0)+\int_0^t\Big(a_{21}^{r}Y_2^{r_2}(s)-a_{11}^{r}Y_1^{r_1}(s)\Big)\,\mrm{d}s+\int_0^t\sqrt{2\alpha Y_1^{r_1}(s)}\,\mrm{d}B(s)\\
&&+\int_0^t\int_0^{Y_1^{r_1}(s-)}\int_{\mbb{M}\setminus A_r}z_1\,\tilde{N}_1(\mrm{d}s,\mrm{d}u,\mrm{d}z)
+\int_0^t\int_0^{Y_2^{r_2}(s-)}\int_{{\mbb{M}_{-1}}\setminus A_r}z_1\,\tilde{N}_2(\mrm{d}s,\mrm{d}u,\mrm{d}z),\\
&Y_2^{r_2}(t)&=Y_2(0)+\int_0^t\Big(b_{11}^{r}Y_1^{r_1}(s)+b_{21}^{r}Y_2^{r_2}(s)\Big)\,\mrm{d}s\\
&&+\int_0^t\int_0^{Y_1^{r_1}(s-)}\int_{\mbb{M}\setminus A_r}z_2\,\tilde{N}_1(\mrm{d}s,\mrm{d}u,\mrm{d}z)
+\int_0^t\int_0^{Y_2^{r_2}(s-)}\int_{{\mbb{M}_{-1}}\setminus A_r}z_2\,\tilde{N}_2(\mrm{d}s,\mrm{d}u,\mrm{d}z).
\eeqnn
Then $\{(Y_1^{r_1}(t),Y_2^{r_2}(t)):t\geq0\}$ is a MSB-process with branching mechanism $(\Phi_1^r,\Phi_2^r)$. It is easy to see that $(Y_1^{r_1}(s),Y_2^{r_2}(s))=(Y_1(s),Y_1(s))$ for $0\leq s<\tau_r$ and
\beqnn
&&\{\tau_r>t\}\\
&&=\Big\{\max\limits_{0<s\leq t}\Delta Y_1(s)\leq r_1,~\max\limits_{0<s\leq t}\Delta Y_2(s)\leq r_2\Big\}\\
&&=\Big\{\int_0^t\int_0^{Y_1(s-)}\int_{A_r}\,N_1(\mrm{d}s,\mrm{d}u,\mrm{d}z)=0,\int_0^t\int_0^{Y_2(s-)}\int_{A_r}\,N_2(\mrm{d}s,\mrm{d}u,\mrm{d}z)=0\Big\}\\
&&=\Big\{\int_0^t\int_0^{Y_1^{r_1}(s-)}\int_{A_r}\,N_1(\mrm{d}s,\mrm{d}u,\mrm{d}z)=0,\int_0^t\int_0^{Y_2^{r_2}(s-)}\int_{A_r}\,N_2(\mrm{d}s,\mrm{d}u,\mrm{d}z)=0\Big\}.\label{5.7}
\eeqnn
Note that $N_1$ and $N_2$ restricted to $(0,\infty)^2\times A_r$ are independent of $\{(Y_1^{r_1}(t),Y_2^{r_2}(t)):t\geq0\}$. It follows that
$$
\mbf{P}_y\{\tau_r>t\}=\mbf{P}_y\,\exp\Big\{-n_1(A_r)\int_0^tY_1^{r_1}(s)\mrm{d}s-n_2(A_r)\int_0^tY_2^{r_2}(s)\mrm{d}s\Big\},
$$
Finally the desired result follows from (\ref{5.1})--(\ref{5.4}).
\qed

\bcorollary\label{co:5.5}
  Suppose that both $n_1$ and $n_2$ have unbounded supports. As $r\rightarrow\mbf{\infty},$ we have
\beqlb
\mbf{P}_y\{\tau_r\leq t\}\sim
\left(\begin{array}{ccc}
y_1,y_2
\end{array}
\right)
\int_0^t \mrm{e}^{(t-s)\mrm{H}}\,\mrm{d}s
\left( \begin{array}{ccc}
n_1(A_r)\\
\\
n_2(A_r)
\end{array}
\right ),\nonumber
\eeqlb
where
\beqlb
\mrm{H}=
\left( \begin{array}{ccc}
-a_{11}&\int_\mbb{M}z_2\,n_1(\mrm{d}z)\\
\\
a_{21}+\int_{\mbb{M}_{-1}}z_1\,n_2(\mrm{d}z)&
\int_{\mbb{M}_{-1}}z_2\,n_2(\mrm{d}z)
\end{array}
\right )\nonumber
\eeqlb
and $\mrm{e}^{(t-s)\mrm{H}}=\sum\limits_{k=0}^{\infty}\frac{(t-s)^k\mrm{H}^k}{k!}$.
\ecorollary

\proof
For $r,q\in[0,\infty)^2,q_i\geq r_i, i=1,2,$ we have obviously $\Phi_i\geq\Phi_i^r\geq\Phi_i^q, i=1,2.$ Then by Proposition 5.2 we see that $\bar{V}_i\geq\bar{V}^{r}_i\geq\bar{V}^{q}_i, i=1,2.$ It follows that
$$
1-\mrm{exp}\Big\{-\langle y,\bar{V}^q(t,n(A_r))\rangle\Big\}\leq\mbf{P}_y\{\tau_r\leq t\}\leq1-\mrm{exp}\Big\{-\langle y,\bar{V}(t,n(A_r))\rangle\Big\}.
$$
Note that $\bar{V}^q(t,0+)=\bar{V}^r(t,0+)=\bar{V}(t,0+)=0.$ Moreover, we can calculate that
\beqlb
&&\frac{\partial}{\partial t}\frac{\partial}{\partial\lambda}\bar{V}_1(t,0)=\mrm{e}^{(1)}-a_{11}\frac{\partial}{\partial\lambda}\bar{V}_1(t,0)
+\int_\mbb{M}z_2\,n_1(\mrm{d}z)\cdot\frac{\partial}{\partial\lambda}\bar{V}_2(t,0),\nonumber\\
&&\frac{\partial}{\partial t}\frac{\partial}{\partial\lambda}\bar{V}_2(t,0)=\mrm{e}^{(2)}+(a_{21}+\int_{\mbb{M}_{-1}}z_1\,n_2(\mrm{d}z))
\frac{\partial}{\partial\lambda}\bar{V}_1(t,0)+\int_{\mbb{M}_{-1}}z_2\,n_2(\mrm{d}z)\cdot\frac{\partial}{\partial\lambda}\bar{V}_2(t,0),\nonumber\\
&&\frac{\partial}{\partial\lambda}\bar{V}_1(0,0)=\frac{\partial}{\partial\lambda}\bar{V}_2(0,0)=0.\nonumber
\eeqlb
We can solve the above equations to get
\beqlb
\left( \begin{array}{ccc}
\frac{\partial}{\partial\lambda}\bar{V}_1(t,0)\\
\\
\frac{\partial}{\partial\lambda}\bar{V}_2(t,0)
\end{array}
\right )
=\int_0^t \mrm{e}^{(t-s)\mrm{H}}\,\mrm{d}s.\label{5.5}
\eeqlb
Similarly we have
\beqlb
\left( \begin{array}{ccc}
\frac{\partial}{\partial\lambda}\bar{V}^q_1(t,0)\\
\\
\frac{\partial}{\partial\lambda}\bar{V}^q_2(t,0)
\end{array}
\right )
=\int_0^t \mrm{e}^{(t-s)\mrm{H}_q}\,\mrm{d}s,\label{5.6}
\eeqlb
where
\beqlb
\mrm{H}_q=
\left( \begin{array}{ccc}
-a^q_{11}&\int_{\mbb{M}\setminus A_q}z_2\,n_1(\mrm{d}z)\\
\\
a_{21}+\int_{{\mbb{M}_{-1}}\setminus A_q}z_1\,n_2(\mrm{d}z)&
\int_{{\mbb{M}_{-1}}\setminus A_q}z_2\,n_2(\mrm{d}z)
\end{array}
\right ).\nonumber
\eeqlb
By (\ref{5.5}) and (\ref{5.6}), as $r\rightarrow\infty,$
\beqnn
&1-\mrm{exp}\Big\{-\langle y,\bar{V}(t,n(A_r))\rangle\Big\}
&\sim\langle y,\bar{V}(t,n(A_r))\rangle\\
&&\sim\left( \begin{array}{ccc}
y_1,y_2
\end{array}
\right)
\int_0^t \mrm{e}^{(t-s)\mrm{H}}\,\mrm{d}s
\left(\begin{array}{ccc}
n_1(A_r)\\
\\
n_2(A_r)
\end{array}
\right )
\eeqnn
and
\beqnn
&1-\mrm{exp}\Big\{-\langle y,\bar{V}^q(t,n(A_r))\rangle\Big\}&\sim\langle y,\bar{V}^q(t,n(A_r))\rangle\nonumber\\
&&\sim \left(\begin{array}{ccc}
y_1,y_2
\end{array}
\right)
\int_0^t \mrm{e}^{(t-s)\mrm{H}_q}\,\mrm{d}s
\left(\begin{array}{ccc}
n_1(A_r)\\
\\
n_2(A_r)
\end{array}
\right).
\eeqnn
Then we complete the proof by noticing $\lim\limits_{q\rightarrow\infty}\mrm{H}_q=\mrm{H}.$
\qed

\bigskip

\section{Exponential ergodicity in Wasserstein distances}

\setcounter{equation}{0}

In order to present our results in this section, we first introduce some notations. Given two probability measures $\mu$ and $\nu$ on $\mbb{M}$, the standard $L^p$-Wasserstein distance $W_p$ for all $p\geq1$ is given by
$$
W_p(\mu,\nu)=\inf_{\Pi\in\mathcal{C}(\mu,\nu)}\bigg(\int_{\mbb{M}\times \mbb{M}}|x-y|^p\,\Pi(\mathrm{d}x,\mathrm{d}y)\bigg)^{1/p},
$$
where $|\cdot|$ denotes the Euclidean norm and $\mathcal{C}(\mu,\nu)$ stands for the set of all coupling measures of $\mu$ and $\nu$, i.e. $\mathcal{C}(\mu,\nu)$ is the collection of measures on $\mbb{M}\times\mbb{M}$ having $\mu$ and $\nu$ as marginals. Denote $\mathcal{P}_{p}(\mbb{M})$ as the set of probability measures having finite moment of order $p$, it is known that
$(\mathcal{P}_{p}(\mbb{M}),W_p)$ becomes a Polish space.

The next theorem gives upper and lower bounds for the variations in the $L^1$-Wasserstein distance $ W_1$ of the transition probabilities of the MSB-process started from two different initial states.

\btheorem\label{th:6.1}
Let $(P_t)_{t\geq0}$ be the transition semigroup defined by (\ref{3.10}). Then for all $x,y\in\mbb{M}$ and $t\geq0$ we have
\beqnn
|\langle x-y,\pi(t,1)\rangle|\leq W_1(\delta_xP_t,\delta_yP_t)\leq\sum_{i=1}^2|x_i-y_i|\pi_i(t,1),
\eeqnn
where $\delta_xP_t(\cdot):=P_t(x,\cdot)$ and $\pi(t,1)$ is defined as in Proposition \ref{pr:3.4} with $\lambda= (1,1)$.
\etheorem

\proof The proof is based on the same idea as that of Theorem 2.2 in \cite{L20}.
By Proposition \ref{pr:3.4}, we see that $\int_{\mbb{M}}(y_1+y_2)\,P_t(x,\mrm{d}y)=\langle x,\pi(t,1)\rangle.$ It follows from Theorem 5.10 in \cite{C04} that
$$
W_1(\delta_xP_t,\delta_yP_t)\geq\int_{\mbb{M}}(z_1+z_2)\,\Big(P_t(x,\mrm{d}z)-P_t(y,\mrm{d}z)\Big)=
\langle x-y,\pi(t,1)\rangle.
$$
Similarly, $W_1(\delta_xP_t,\delta_yP_t)\geq\langle y-x,\pi(t,1)\rangle.$ Then the first inequality follows. On the other hand, for $x,y\in\mbb{M},$ let $(x-y)_{\pm}:=((x_1-y_1)_{\pm},(x_2-y_2)_{\pm})$, and $x\wedge y:=x-(x-y)_{+}=y-(x-y)_{-}.$
Let $P_t(x,y,\mrm{d}\eta_1,\mrm{d}\eta_2)$ be the image of the product measure
$$
P_t(x\wedge y,\mrm{d}\gamma_0)P_t((x-y)_{+},\mrm{d}\gamma_1)P_t((x-y)_{-},\mrm{d}\gamma_2)
$$
under the mapping $(\gamma_0,\gamma_1,\gamma_2)\mapsto(\eta_1,\eta_2):=(\gamma_0+\gamma_1,\gamma_0+\gamma_2)$. It's not hard to see that
$P_t(x,y,\mrm{d}\eta_1,\mrm{d}\eta_2)$ is a coupling of $P_t(x,\mrm{d}\eta_1)$ and $P_t(y,\mrm{d}\eta_2)$. Then
\beqnn
&&W_1(\delta_xP_t,\delta_y P_t)
\leq\int_{\mbb{M}^2}|\eta_1-\eta_2|\,P_t(x,y,\mrm{d}\eta_1,\mrm{d}\eta_2)\\
&&~~~~~~~~~~~~~~~~~~~
\leq\int_{\mbb{M}}P_t((x-y)_{+},\mrm{d}\gamma_1)\int_{\mbb{M}}(\gamma_{11}+\gamma_{12}+\gamma_{21}+\gamma_{22})\,P_t((x-y)_{-},\mrm{d}\gamma_2)\\
&&~~~~~~~~~~~~~~~~~~~
=\int_{\mbb{M}}(\zeta_{1}+\zeta_{2})\,P_t((|x_1-y_1|,|x_2-y_2|),\mrm{d}\zeta)\\
&&~~~~~~~~~~~~~~~~~~~
=\sum_{i=1}^2|x_i-y_i|\pi_i(t,1),
\eeqnn
where we have used the branching property $P_t(a,\cdot)\ast P_t(b,\cdot)=P_t(a+b,\cdot)$
for all $a,b\in\mbb{M},t\geq0$ in the third row. Therefore the proof is finished.
 \qed

Based on Theorem \ref{th:6.1}, we can establish the exponential ergodicity with respect to $W_1$. Recalling that a $2\times2$ matrix $H=[H_{ij}]_{2\times2}$ in Corollary \ref{co:5.5} is defined as follows:
\beqlb
\mrm{H}=
\left( \begin{array}{ccc}
-a_{11}&\int_\mbb{M}z_2\,n_1(\mrm{d}z)\\
\\
a_{21}+\int_{\mbb{M}_{-1}}z_1\,n_2(\mrm{d}z)&
\int_{\mbb{M}_{-1}}z_2\,n_2(\mrm{d}z)
\end{array}
\right ),\nonumber
\eeqlb
we have the following result:

\btheorem\label{th:6.2}
Assume that $H_{11}H_{22}-H_{12}H_{21}>0$ and $H_{11}+H_{22}<0.$ Then there exist $\lambda,\vartheta>0$ such that for any $t\geq0$ and $x,y\in\mbb{M}$,
$$
W_1(\delta_xP_t,\delta_yP_t)\leq\vartheta|x-y|\mrm{e}^{-\lambda t}.
$$
\etheorem

\proof
By assumption, it is easy to see that
$$
\lambda^2-(H_{11}+H_{22})\lambda+H_{11}H_{22}-H_{12}H_{21}=0
$$
has two different roots: $\lambda_1=2^{-1}(H_{11}+H_{22}+\sqrt{\Delta}),$
$\lambda_2=\lambda_1-\sqrt{\Delta}$ and $\lambda_2<\lambda_1<0$, where $\Delta=(H_{11}-H_{22})^2+4H_{12}H_{21}>0$. If $H_{12}=H_{21}=0$, then $\lambda_1=H_{11}$ and $\lambda_2=H_{22}$. By Proposition \ref{pr:3.4}, $\pi_i(t,1)=\mrm{e}^{\lambda_it}$ for $i=1,2$, and the desired result follows.
Next we only consider $H_{12}>0$. We can calculate that
\beqnn
&&\pi_1(t,1)=\frac{H_{11}+H_{12}-\lambda_2}{\sqrt{\Delta}}\mrm{e}^{\lambda_1t}
+\frac{\lambda_1-H_{11}-H_{12}}{\sqrt{\Delta}}\mrm{e}^{\lambda_2t}\\
&&~~~~~~~~:=\theta_{11}\mrm{e}^{\lambda_1t}+\theta_{12}\mrm{e}^{\lambda_2t},\\
&&\pi_2(t,1)=\frac{(H_{11}+H_{12}-\lambda_2)(\lambda_1-H_{11})}{\sqrt{\Delta}H_{12}}\mrm{e}^{\lambda_1t}
+\frac{(\lambda_1-H_{11}-H_{12})(\lambda_2-H_{11})}{\sqrt{\Delta}H_{12}}\mrm{e}^{\lambda_2t}\\
&&~~~~~~~~:=\theta_{21}\mrm{e}^{\lambda_1t}+\theta_{22}\mrm{e}^{\lambda_2t}.
\eeqnn
It is easy to see that $\theta_{11},\theta_{21}>0$. It follows from Theorem \ref{th:6.1} that
\beqnn
&&W_1(\delta_xP_t,\delta_yP_t)\leq|x_1-y_1|(\theta_{11}\mrm{e}^{\lambda_1t}+\theta_{12}\mrm{e}^{\lambda_2t})
+|x_2-y_2|(\theta_{21}\mrm{e}^{\lambda_1t}+\theta_{22}\mrm{e}^{\lambda_2t})\\
&&~~~~~~~~~~~~~~~~~~~
\leq|x-y|(|\theta_{11}|+|\theta_{21}|)\mrm{e}^{\lambda_1t}+|x-y|(|\theta_{12}|+|\theta_{22}|)\mrm{e}^{\lambda_2t}\\
&&~~~~~~~~~~~~~~~~~~~
\leq(\theta_{11}+\theta_{21}+|\theta_{12}|+|\theta_{22}|)|x-y|\mrm{e}^{\lambda_1t},
\eeqnn
we obtain the desired result by setting $\vartheta=\theta_{11}+\theta_{21}+|\theta_{12}|+|\theta_{22}|>0$ and $\lambda=-\lambda_1>0$.
\qed

\bcorollary\label{co:6.3}
Assume the conditions of Theorem \ref{th:6.2} hold. Then there exist a unique $\pi\in\mathcal{P}_{1}(\mbb{M})$ and $\vartheta,\lambda>0$ such that for any $x\in\mbb{M}$ and $t\geq0$
$$
W_1(\delta_xP_t, \pi)\leq\vartheta W_1(\delta_x, \pi)\mrm{e}^{-\lambda t}.
$$
\ecorollary

\proof
By Theorem \ref{th:7.5} below, there exists a unique invariant measure. Arguing similarly to the proof of Theorem 3.2 in \cite{FJKR19}, one can see that $\pi\in\mathcal{P}_{1}(\mbb{M}),$  and the desired assertion is easily obtained by Theorem \ref{th:6.2}.
\qed

\section{MSBI-processes}
\setcounter{equation}{0}
Suppose that $\Phi_1, \Phi_2$ are two functions on $[0,\infty)^2$ defined as in (\ref{3.3})--(\ref{3.4}), and there exists function $\Psi$ on $[0,\infty)^2$ defined by:
\beqlb\label{7.1}
\Psi(\lambda_1,\lambda_2)=b\lambda_1+\int_{\mbb{M}}(1-\mrm{e}^{-\langle\lambda,z\rangle})m(\mrm{d}z),\quad \lambda\in\mbb{R}^2_{+},
\eeqlb
where $b>0$ and $m$ is a $\sigma-$finite measure on $\mbb{M}$ supported by $\mbb{M}\setminus\{\mathbf{0}\}$ such that
$$
\int_{\mbb{M}}(1\wedge z_1+1\wedge z_2)\,m(\mrm{d}z)<\infty.
$$

A Markov process $\{Z(t)=(Z_1(t),Z_2(t)):t\geq0\}$ is called a MSBI-process on $\mbb{M},$ if it has transition semigroup $(P^\gamma_t)_{t\geq0}$ uniquely determined by:
\beqlb
\int_\mbb{M}\mrm{e}^{-\langle\lambda,y\rangle}P^\gamma_t(x,\mrm{d}y)=\exp\Big\{-\langle x,V(t,\lambda)\rangle-\int_0^t\Psi(V(s,\lambda))\,\mrm{d}s\Big\},\quad x\in \mbb{M},\lambda\in\mbb{R}^2_{+},\label{7.2}
\eeqlb
where $V(t,\lambda)=(V_1(t,\lambda),V_2(t,\lambda))$ takes values on $\mbb{R}_{+}^2$ and satisfies (\ref{3.8}). One can see that the semigroup defined by (\ref{7.2}) is a Feller semigroup, then the MSBI-process has a $\mrm{c\grave{a}dl\grave{a}g}$ realization. We can also establish the similar result of Theorem \ref{th:6.1} for MSBI-processes, indeed, we have the following:

\btheorem\label{th:7.1}
Let $(P^\gamma_t)_{t\geq0}$ be the transition semigroup defined by (\ref{7.2}). Assume that $\int_{\mbb{M}}(z_1+z_2)\,m(\mrm{d}z)<\infty.$ Then for $t\geq0$ and $x,y\in\mbb{M}$ we have
\beqnn
|\langle x-y,\pi(t,1)\rangle|\leq W_1(\delta_xP^\gamma_t,\delta_yP^\gamma_t)\leq\sum_{i=1}^2|x_i-y_i|\pi_i(t,1),
\eeqnn
where $\pi(t,1)$ is defined as in Proposition \ref{pr:3.4} with $\lambda= (1,1)$.
\etheorem

\smallskip
\proof The proof is based on the same idea as that of Theorem 4.1 in \cite{L20}. One can see that
\beqnn
&&\int_{\mbb{M}}(y_1+y_2)\,P^\gamma_t(x,\mrm{d}y)=\langle x,\pi(t,1)\rangle+(b+\int_{\mbb{M}}z_1\,m(\mrm{d}z))\int_0^t\pi_1(s,1)\,\mrm{d}s\\
&&~~~~~~~~~~~~~~~~~~~~~~~~~~~~~~~~
+\int_{\mbb{M}}z_2\,m(\mrm{d}z)\int_0^t\pi_2(s,1)\,\mrm{d}s,
\eeqnn
which yields that
$$
W_1(\delta_xP^\gamma_t,\delta_yP^\gamma_t)\geq\int_{\mbb{M}}(z_1+z_2)\,(P^\gamma_t(x,\mrm{d}z)-P^\gamma_t(y,\mrm{d}z))
=\langle x-y,\pi(t,1)\rangle,
$$
similarly, $W_1(\delta_xP^\gamma_t,\delta_yP^\gamma_t)\geq\langle y-x,\pi(t,1)\rangle$, and the first inequality follows. Next, we want to construct a coupling measure of $P^\gamma_t(x,\cdot)$ and $P^\gamma_t(y,\cdot)$. It is known that there exists a family of probability measures $(\gamma_t)_{t\geq0}$
such that $P^\gamma_t(x,\cdot)=P_t(x,\cdot)\ast\gamma_t(\cdot)$ for $t\geq0,x\in\mbb{M}$, and
$$
\int_{\mbb{M}}\mrm{e}^{-\langle\lambda,y\rangle}\,\gamma_t(\mrm{d}y)=\exp\Big\{-\int_0^t\Psi(V(s,\lambda))\,\mrm{d}s\Big\},
$$
we call $(\gamma_t)_{t\geq0}$ a skew convolution semigroup associated with $(P_t)_{t\geq0}$; see, e.g., Chapter 9 in \cite{L11}. Let $P_t(x,y,\mrm{d}\eta_1,\mrm{d}\eta_2)$ be the coupling measure of $P_t(x,\mrm{d}\eta_1)$ and $P_t(y,\mrm{d}\eta_2)$ constructed in the proof of Theorem \ref{th:6.1} and $P^\gamma_t(x,y,\mrm{d}\sigma_1,\mrm{d}\sigma_2)$ be the image of $\gamma_t(\mrm{d}\eta_0)P_t(x,y,\mrm{d}\eta_1,\mrm{d}\eta_2)$ under the mapping $(\eta_0,\eta_1,\eta_2)\mapsto(\sigma_1,\sigma_2)=(\eta_0+\eta_1,\eta_0+\eta_2)$. By the relation
$P^\gamma_t(x,\cdot)=P_t(x,\cdot)\ast\gamma_t(\cdot)$ we see that $P^\gamma_t(x,y,\mrm{d}\sigma_1,\mrm{d}\sigma_2)$ is a coupling measure of $P^\gamma_t(x,\mrm{d}\sigma_1)$ and $P^\gamma_t(y,\mrm{d}\sigma_2)$. It follows that
\beqnn
&&W_1(\delta_xP^\gamma_t,\delta_yP^\gamma_t)\leq\int_{\mbb{M}^2}|\sigma_1-\sigma_2|\,P^\gamma_t(x,y,\mrm{d}\sigma_1,\mrm{d}\sigma_2)\\
&&~~~~~~~~~~~~~~~~~~~~~
=\int_{\mbb{M}}\,\gamma_t(\mrm{d}\eta_0)\int_{\mbb{M}^2}|\eta_1-\eta_2|\,P_t(x,y,\mrm{d}\eta_1,\mrm{d}\eta_2)\\
&&~~~~~~~~~~~~~~~~~~~~~
=\int_{\mbb{M}^2}|\eta_1-\eta_2|\,P_t(x,y,\mrm{d}\eta_1,\mrm{d}\eta_2)
\leq\sum_{i=1}^2|x_i-y_i|\pi_i(t,1),
\eeqnn
where the last inequality follows from Theorem \ref{th:6.1}.
 \qed

By a similar argument of Theorem \ref{th:6.2}, we have the following:

\btheorem\label{th:7.2}

Assume that $H_{11}H_{22}-H_{12}H_{21}>0$ and $H_{11}+H_{22}<0.$ Then there exist $\lambda,\vartheta>0$ such that for any $t\geq0$ and $x,y\in\mbb{M}$,
$$
W_1(\delta_xP^{\gamma}_t,\delta_yP^{\gamma}_t)\leq\vartheta|x-y|\mrm{e}^{-\lambda t}.
$$
\etheorem

\subsection{The construction of MSBI-processes by stochastic equations}

We now give a construction of MSBI-processes by stochastic equations. Let us consider the following stochastic equation system:
\beqlb
&&Z_1(t)=Z_1(0)+\int_0^t\Big(b-a_{11}Z_1(s)+a_{21}Z_2(s)\Big)\,\mrm{d}s
+\int_0^t\sqrt{2\alpha Z_1(s)}\,\mrm{d}B(s)\nonumber \\
&&\qquad~~~~~~~~~~~
+\int_0^t\int_0^{Z_1(s-)}\int_\mbb{M}z_1\,\tilde{N}_1(\mrm{d}s,\mrm{d}u,\mrm{d}z)
+\int_0^t\int_\mbb{M}z_1\,M(\mrm{d}s,\mrm{d}z)\nonumber\\
&&\qquad~~~~~~~~~~~
+\int_0^t\int_0^{Z_2(s-)}\int_{\mbb{M}_{-1}}z_1\,N_2(\mrm{d}s,\mrm{d}u,\mrm{d}z),
\label{7.3}\\
&&\qquad \nonumber\\
&&Z_2(t)=Z_2(0)+\int_0^t\int_0^{Z_1(s-)}\int_\mbb{M}z_2\,N_1(\mrm{d}s,\mrm{d}u,\mrm{d}z)
+\int_0^t\int_\mbb{M}z_2\,M(\mrm{d}s,\mrm{d}z)\nonumber\\
&&\qquad~~~~~~~~~~~
+\int_0^t\int_0^{Z_2(s-)}\int_{\mbb{M}_{-1}}z_2\,N_2(\mrm{d}s,\mrm{d}u,\mrm{d}z),
\label{7.4}
\eeqlb
where $b\geq0,M(\mrm{d}s,\mrm{d}z)$ is a Poisson random measure on $[0,\infty)\times\mbb{M}$ with intensity measure $\mrm{d}sm(\mrm{d}z),$ other coefficients are the same in section 4. Furthermore, we assume that those random elements are independent of each other. By a modification of the proof of section 4 as well as in \cite{M13}, we see that (\ref{7.3})--(\ref{7.4}) has a unique strong solution and it is a MSBI-process with branching mechanism $(\Phi_1,\Phi_2)$ defined by (\ref{3.3})--(\ref{3.4}) and immigration mechanism $\Psi$ defined by (\ref{7.1}).

\subsection{Stationary distribution}
In order to characterize the stationary distribution of MSBI-processes, we need to estimate the upper bound and lower bound of $|V(t,\lambda)|$ for $t>0,\lambda\in\mbb{R}_{+}^2$, which will play an important role in the sequel.

\begin{lemma}\label{le:7.3}
Let $(Y_t)_{t\geq0}$ be a MSB-process with semigroup $(P_t)_{t\geq0}$ satisfying (\ref{3.10}). Let $H=[H_{ij}]_{2\times2}$ be a $2\times2$ matrix defined as in Corollary \ref{co:5.5}. Suppose that all the eigenvalues of $H$ have strictly negative real parts. Then there exist some strictly positive constants $c_1(\lambda)$ and $c_2$, where $c_1$ depends on $\lambda$ such that
$$
|V(t,\lambda)|\leq c_1(\lambda)\exp\{-c_2t\},\quad \lambda\in\mbb{R}_{+}^2,\quad t\geq0.
$$
\end{lemma}

 \proof
We follow the same calculations in Proposition \ref{pr:3.4} to see that
\beqlb
\left( \begin{array}{ccc}
\frac{\partial V_1(t,0+)}{\partial\lambda_1}\\
\\
\frac{\partial V_2(t,0+)}{\partial\lambda_1}
\end{array}
\right )
=\mrm{e}^{tH}
\left( \begin{array}{ccc}
1\\
\\
0
\end{array}
\right ),\nonumber
\eeqlb
and so
\beqlb
\int_{\mbb{M}}y_1\,P_t(x,\mrm{d}y)=
\left( \begin{array}{ccc}
x_1\\
\\
x_2
\end{array}
\right )^{\mrm{T}}
\mrm{e}^{tH}
\left( \begin{array}{ccc}
1\\
\\
0
\end{array}
\right ).\nonumber
\eeqlb
Similarly,
\beqlb
\int_{\mbb{M}}y_2\,P_t(x,\mrm{d}y)=
\left( \begin{array}{ccc}
x_1\\
\\
x_2
\end{array}
\right )^{\mrm{T}}
\mrm{e}^{tH}
\left( \begin{array}{ccc}
0\\
\\
1
\end{array}
\right ).\nonumber
\eeqlb
By Jensen's inequality we deduce that for all $x=(x_1,x_2)\in\mbb{M}$
\beqlb
\left( \begin{array}{ccc}
x_1\\
\\
x_2
\end{array}
\right )^{\mrm{T}}\left( \begin{array}{ccc}
V_1(t,\lambda)\\
\\
V_2(t,\lambda)
\end{array}
\right )\leq
\left( \begin{array}{ccc}
x_1\\
\\
x_2
\end{array}
\right )^{\mrm{T}}
\mrm{e}^{tH}
\left( \begin{array}{ccc}
\lambda_1\\
\\
\lambda_2
\end{array}
\right ).\nonumber
\eeqlb
Since all the eigenvalues of $H$ have strictly negative real parts, there exist some strictly positive $c,c_2>0$ such that for all $t>0$
$$
\|\mrm{e}^{tH}\|:=\sup_{|x|=1}|\mrm{e}^{tH}x|\leq c\mrm{e}^{-c_2t};
$$
see,e.g., equation (2.8) in \cite{SY84}, which implies that $|V(t,\lambda)|\leq|\lambda|c\mrm{e}^{-c_2t},$ we finish the proof by setting $c_1(\lambda)=|\lambda|c.$
 \qed

\begin{lemma}\label{le:7.4}
Under the conditions of Lemma \ref{le:7.3}, for every $\lambda\in\mbb{R}_{+}^2$, there exist two strictly positive constants $A(\lambda)$ and $B(\lambda)$ such that
$$
V_1(t,\lambda)\geq\lambda_1\mrm{e}^{-A(\lambda)t},\quad
V_2(t,\lambda)\geq\lambda_2\mrm{e}^{-B(\lambda)t},\qquad t\geq0.
$$
\end{lemma}

\proof
\beqnn
&&V_1(t,\lambda)=\lambda_1+\int_0^t\Big\{H_{11}V_1(s,\lambda)+H_{12}V_2(s,\lambda)-\alpha V_1^2(s,\lambda)\\
&&~~~~~~~~~~~~~~~~~~~~~~~
-\int_{\mbb{M}}(\mrm{e}^{-\langle V(s,\lambda),z\rangle}-1+\langle V(s,\lambda),z\rangle)\,n_1(\mrm{d}z)\Big\}\,\mrm{d}s,\\
&&V_2(t,\lambda)=\lambda_2+\int_0^t\Big\{a_{21}V_1(s,\lambda)
+\int_{\mbb{M}_{-1}}(1-\mrm{e}^{-\langle V(s,\lambda),z\rangle})\,n_2(\mrm{d}z)\Big\}\,\mrm{d}s.
\eeqnn
It follows from Lemma \ref{le:7.3}, the comparison theorem and the fact
$$
\mrm{e}^{-\lambda x}-1+\lambda x\leq\Big(\frac{\lambda^2}{2}+\lambda\Big)(x\wedge x^2),\quad x,\lambda\geq0
$$
that there exists $A(\lambda)=|H_{11}-\kappa-(\alpha+\kappa/2)c_1(\lambda)|$, where $\kappa=\int_{\mbb{M}}z_1\wedge z_1^2\,n_1(\mrm{d}z)$ such that
$$
V_1(t,\lambda)\geq\lambda_1\mrm{e}^{-A(\lambda)t},\quad t\geq0,\quad \lambda\in\mbb{R}_{+}^2.
$$
On the other hand, let $\theta=n_2(\mbb{R}_{+}\times\{-1\})<\infty,$
\beqlb
&&V_2(t,\lambda)\geq\lambda_2-\theta\int_0^t(\mrm{e}^{V_2(s,\lambda)}-1)\,\mrm{d}s\nonumber\\
&&~~~~~~~~~~
\geq\lambda_2-2\theta\int_0^t\Big(\mrm{e}^{c_1(\lambda)\mrm{e}^{-c_2s}}V_2(s,\lambda)\Big)\,\mrm{d}s\nonumber\\
&&~~~~~~~~~~
\geq\lambda_2-2\theta\mrm{e}^{c_1(\lambda)}\int_0^tV_2(s,\lambda)\,\mrm{d}s,\nonumber
\eeqlb
by the comparison theorem we deduce that $V_2(t,\lambda)\geq\lambda_2\mrm{e}^{-2\theta\mrm{e}^{c_1(\lambda)}t},$ and we obtain the desired result by setting $B(\lambda)=2\theta\mrm{e}^{c_1(\lambda)}.$
\qed

We now give our main result.

\btheorem\label{th:7.5}
Let $(Z_t)_{t\geq0}$ be a MSBI-process with semigroup $(P^\gamma_t)_{t\geq0}$ satisfying (\ref{7.2}). Suppose that all the eigenvalues of $H$ have strictly negative real parts. Then $P^\gamma_t(x,\cdot)$ converges to a probability measure $\pi$ on $\mbb{M}$ as $t\rightarrow\infty$ for all $x\in\mbb{M}$ if and only if
\beqlb
\int_{\{|z|\geq1\}}\log|z|\,m(\mrm{d}z)<\infty.
\eeqlb
\etheorem

\proof
By Lemma \ref{le:7.3} we have $|V(t,\lambda)|\rightarrow0$ as $t\rightarrow\infty.$ Suppose that $(Z_t)_{t\geq0}$ has a stationary distribution $\pi$, one can see that
$$
\int_{\mbb{M}}\mrm{e}^{-\langle\lambda,y\rangle}\,\pi(\mrm{d}y)
=\exp\Big\{-\int_0^{\infty}\Psi(V(s,\lambda))\,\mrm{d}s\Big\},\quad \lambda\in\mbb{R}_{+}^2,
$$
which implies that $\int_0^{\infty}\Psi(V(s,\lambda))\,\mrm{d}s<\infty$ for all $\lambda\in\mbb{R}_{+}^2.$
And so
$$
\int_0^{\infty}\mrm{d}s\int_{\{|z|\geq1\}}\Big(1-\mrm{e}^{-(\lambda_1\wedge\lambda_2)\exp\{-[A(\lambda)+ B(\lambda)]s\}(z_1+z_2)}\Big)\,m(\mrm{d}z)<\infty.
$$
If we set $C(\lambda):=|H_{11}|+\vartheta+(\alpha+\vartheta/2)c_1(\lambda)+2\theta\mrm{e}^{c_1(\lambda)}>0$
for all $\lambda\in\mbb{R}_{+}^2,$ then $C(\lambda)\geq A(\lambda)+B(\lambda)$. Choose a proper $\tilde{\lambda}$ such that $\tilde{\lambda}_1\wedge\tilde{\lambda}_2>0$ and let $t:=\mrm{e}^{-C(\tilde{\lambda})s}|z|,$ we have $\mrm{d}t=-C(\tilde{\lambda})t\,\mrm{d}s,$ and
$$
\int_{\{|z|\geq1\}}m(\mrm{d}z)\int_0^{|z|}\frac{1-\mrm{e}^{-(\tilde{\lambda}_1\wedge\tilde{\lambda}_2)t}}{t}\,\mrm{d}t<\infty,
$$
which yields that
$$
\int_{\{|z|\geq1\}}\log|z|\,m(\mrm{d}z)<\infty
$$
since
$$
\int_0^{|z|}\frac{1-\mrm{e}^{-(\tilde{\lambda}_1\wedge\tilde{\lambda}_2)t}}{t}\,\mrm{d}t\sim\log|z| \quad \mrm{as}~|z|\rightarrow\infty.
$$
On the other hand, it suffices to prove $\int_0^{\infty}\int_{\mbb{M}}(1-\mrm{e}^{-\langle V(s,\lambda),z\rangle})\,\mrm{d}s\,m(\mrm{d}z)<\infty$ for all $\lambda\in{R}_{+}^2$ provided $\int_{\{|z|\geq1\}}\log|z|\,m(\mrm{d}z)<\infty.$ From Lemma \ref{le:7.3} and Fubini's theorem
\beqnn
&&\int_0^{\infty}\,\mrm{d}s\int_{\mbb{M}}\Big(1-\mrm{e}^{-\langle V(s,\lambda),z\rangle}\Big)\,m(\mrm{d}z)\leq\int_0^{\infty}\,\mrm{d}s\int_{\mbb{M}}\Big(1-\mrm{e}^{-c_1(\lambda)\mrm{e}^{-c_2s}(z_1+z_2)}\Big)\,m(\mrm{d}z)\\
&&~~~~~~~~~~~~~~~~~~~~~~~~~~~~~~~~~~~~~~~~~~~~~~
=\int_0^{\infty}\,\mrm{d}s\int_{\{|z|<1\}}\Big(1-\mrm{e}^{-c_1(\lambda)\mrm{e}^{-c_2s}(z_1+z_2)}\Big)\,m(\mrm{d}z)\\
&&~~~~~~~~~~~~~~~~~~~~~~~~~~~~~~~~~~~~~~~~~~~~~~
+\int_0^{\infty}\,\mrm{d}s\int_{\{|z|\geq1\}}\Big(1-\mrm{e}^{-c_1(\lambda)\mrm{e}^{-c_2s}(z_1+z_2)}\Big)\,m(\mrm{d}z) \\
&&~~~~~~~~~~~~~~~~~~~~~~~~~~~~~~~~~~~~~~~~~~~~~~
:=I_{*}(\lambda)+I^{*}(\lambda).
\eeqnn
For $I_{*}(\lambda),$ by a change of variables $t:=c_1(\lambda)\mrm{e}^{-c_2s}(z_1+z_2)$ we get
\beqnn
&&I_{*}(\lambda)=\frac{1}{c_2}\int_{\{|z|<1\}}\,m(\mrm{d}z)\int_0^{c_1(\lambda)(z_1+z_2)}\frac{1-\mrm{e}^{-t}}{t}\,\mrm{d}t\\
&&~~~~~~~
\leq\frac{c_1(\lambda)}{c_2}\int_{\{|z|<1\}}(z_1+z_2)\,m(\mrm{d}z)<\infty,
\eeqnn
where the last inequality follows from
$\int_{\mbb{M}}(1\wedge z_1+1\wedge z_2)\,m(\mrm{d}z)<\infty.$
Moreover, by a change of variables $t:=c_1(\lambda)\mrm{e}^{-c_2s}|z|,$ we have
\beqlb
&&I^{*}(\lambda)=\frac{1}{c_2}\int_{\{|z|\geq1\}}\,m(\mrm{d}z)
\int_{0}^{c_1(\lambda)|z|}\frac{1-\mrm{e}^{-t\frac{(z_1+z_2)}{|z|}}}{t}\,\mrm{d}t\nonumber\\
&&~~~~~~~
\leq\frac{1}{c_2}\int_{\{|z|\geq1\}}\,m(\mrm{d}z)
\int_{0}^{c_1(\lambda)|z|}\frac{1-\mrm{e}^{-2t}}{t}\,\mrm{d}t,\nonumber
\eeqlb
which implies that $I^{*}(\lambda)<\infty$ by observing
$$
\int_0^{c_1(\lambda)|z|}\frac{1-\mrm{e}^{-2t}}{t}\,\mrm{d}t~~\sim~~\log|z|,\quad |z|\rightarrow\infty
$$
and $\int_{\{|z|\geq1\}}\log|z|\,m(\mrm{d}z)<\infty$.
\qed

\bcorollary\label{co:7.6}
Assume $H_{11}H_{22}-H_{12}H_{21}>0$ and $H_{11}+H_{22}<0$ hold. Moreover, suppose that L\'evy measure $m$ satisfies $\int_{\{|z|>1\}}|z|\,m(\mrm{d}z)<\infty$. Then there exist $\lambda,\vartheta>0$ and unique $\pi\in\mathcal{P}_{1}(\mbb{M})$ such that for any $t\geq0$ and $x\in\mbb{M}$,
$$
W_1(\delta_xP^{\gamma}_t,\pi)\leq\vartheta W_1(\delta_x,\pi)\mrm{e}^{-\lambda t}.
$$
\ecorollary

\proof It follows from Theorem \ref{th:7.5} and assumptions that there exists a unique stationary distribution $\pi$. We can derive easily that $\mbf{E}_x[|Z_t|]<\infty$ for all $t\geq0$ and $x\in\mbb{M}$ by the assumption $\int_{\{|z|>1\}}|z|\,m(\mrm{d}z)<\infty$. By a modification of the proof of Corollary \ref{co:6.3}, we have $\pi\in\mathcal{P}_{1}(\mbb{M})$ and the desired result follows from Theorem \ref{th:7.2}.
\qed


\end{document}